\documentclass[final,onefignum,onetabnum]{siamart171218}

\usepackage{tikz}
\usetikzlibrary{cd}
\usepackage{lipsum}
\usepackage{amsfonts}
\usepackage{graphicx}
\usepackage{epstopdf}
\usepackage{algorithmic}
\ifpdf
  \DeclareGraphicsExtensions{.eps,.pdf,.png,.jpg}
\else
  \DeclareGraphicsExtensions{.eps}
\fi

\usepackage[top=3.375cm,
            bottom=3.375cm,
            inner=4.275cm,
            outer=4.275cm]{geometry}

\newsiamremark{remark}{Remark}
\newsiamremark{hypothesis}{Hypothesis}
\crefname{hypothesis}{Hypothesis}{Hypotheses}
\newsiamthm{claim}{Claim}
\newsiamthm{assumption}{Assumption}
\newsiamremark{example}{Example}
\newsiamthm{conjecture}{Conjecture}

\headers{Mapping Holomorphic Sobolev Spaces}{D.S. Winterrose}

\title{Algebras of Pseudo-Differential Operators \\ Acting on Holomorphic Sobolev Spaces}

\author{David Scott Winterrose \thanks{Department of Applied Mathematics and Computer Science, Technical University of Denmark, Kgs. Lyngby, Denmark (\email{dawin@dtu.dk}).}}

\usepackage{amsopn}
\usepackage{mathrsfs}

\ifpdf
\hypersetup{
  pdftitle={Algebras of Pseudo-Differential Operators Acting on Holomorphic Sobolev Spaces},
  pdfauthor={David Winterrose}
}
\fi

\newcommand{\Op}{\textnormal{Op}}

\newcommand{\Diff}{\textnormal{Diff}}
\newcommand{\PsiPHG}{\Psi_\textnormal{phg}}

\newcommand{\Real}{\textnormal{Re}\,}
\newcommand{\Imag}{\textnormal{Im}\,}

\newcommand{\Tr}{\textnormal{Tr}}

\newcommand{\Ad}{\textnormal{Ad}}

\begin{document}

\maketitle

\begin{abstract}
We search for pseudo-differential operators acting on holomorphic Sobolev spaces. 
The operators should mirror the standard Sobolev mapping property in the holomorphic analogues.
The setting is a closed real-analytic Riemannian manifold, or Lie group with a bi-invariant metric, and the holomorphic Sobolev spaces are defined on a fixed Grauert tube about the core manifold.
We find that every pseudo-differential operator in the commutant of the Laplacian is of this kind.
Moreover, so are all the operators in the commutant of certain analytic
pseudo-differential operators, but for more general tubes, provided that 
an old statement of Boutet de Monvel holds true generally.
In the Lie group setting, we find even larger algebras, and characterize all their elliptic elements.
These latter algebras are determined by global matrix-valued symbols.
\end{abstract}



\section{Introduction}
Let $(M,g)$ be a compact real-analytic Riemannian manifold.
This manifold always admits a (non-unique) Bruhat-Whitney complexification $M_\mathbb{C}$, and, inside $M_\mathbb{C}$, the family of Grauert tubes $\{M_\epsilon\}_{\epsilon < \epsilon_0}$. See e.g. \cite{LeichtnamGolseStenzel2015} for further details.
Let us denote by $\Psi(M)$ the algebra of standard pseudo-differential operators on $M$, and by $\mathcal{O}(M_\epsilon)$ the holomorphic functions on $M_\epsilon$ ($\mathcal{O}(N)$ on a complex manifold $N$).
Then, if $P \in \Psi(M)$, we
seek a commutative diagram
\[
\begin{tikzcd}
H_1 \arrow[r, "\widetilde{P}"] \arrow[d, "\mathcal{R}_\epsilon"] & H_2 \arrow[d, "\mathcal{R}_\epsilon"]  \\
 C^\infty(M)  \arrow[r, "P"] & C^\infty(M)   
\end{tikzcd}
\]
where $H_1,H_2 \subset \mathcal{O}(M_\epsilon)$ are function spaces, $\mathcal{R}_\epsilon$ is the restriction map
\begin{align*}
\mathcal{R}_\epsilon : \mathcal{O}(M_\epsilon) \to C^\omega(M) : f \mapsto f|_M.
\end{align*}
That is, we seek necessary and/or sufficient conditions for the existence of such a $\widetilde{P}$, and we would also like the two function spaces $H_1$ and $H_2$ above to be Hilbert spaces.
An operator $P$ with this property will then be called $\epsilon$-holomorphically extendible.
Our question will lead us to commutative diagrams of the form
\[
\begin{tikzcd}[row sep=large, column sep = huge]
HH^s(M_\epsilon) \arrow[r, "\widetilde{P}"] \arrow[d, "\mathcal{R}_\epsilon"] & HH^{s-d}(M_\epsilon) \arrow[d, "\mathcal{R}_\epsilon"]  \\
H^{s}(M) \arrow[r, "P"] & H^{s-d}(M)   
\end{tikzcd}
\]
where $HH^s(M_\epsilon)$ are the holomorphic Sobolev spaces indexed by $s\in \mathbb{R}$, defined later,
and $P\in \Psi^d(M)$ is a (standard) pseudo-differential operator of degree $d\in \mathbb{R}$.

\newpage
Why do we ask this? We seek a parallel to the Boutet de Monvel calculus \cite{GGPsiDBO, BoutetDeMonvel1971}, but adapted to solve real-analytic elliptic boundary value problems in Hilbert spaces of holomorphic functions. This research path was inspired by Karamehmedovi\'c \cite{Karamehmedovic2015}, and further additions in Winterrose \cite{Winterrose2021}. Unfortunately, no such calculus is known. It is not clear if it even exists. As a beginning, we look for $P$ with the above property, and, especially, parametrices with a similar extension property if $P$ is also elliptic.
 The $\epsilon$-holomorphic extension property can be viewed as an abstract global analogue of the property in Winterrose \cite{WinterroseArticle2021} for (local) analytic pseudo-differential operators.  
This work, and the reason for undertaking it, is further detailed in \cite{Winterrose2021}.\\

Consider a compact Lie group $G$ with Lie algebra $\mathfrak{g}$ of left-invariant vector fields. 
It admits a complexification $G_\mathbb{C}$ itself a complex Lie group \cite[pp. 110, Lemma 3]{Hall1994}, and the natural Cartan decomposition diffeomorphism
\begin{align*}
    G \times \mathfrak{g} \to G_\mathbb{C} : (x,Y) \mapsto \exp(iY)x,
\end{align*}
where $i$ comes from the complex structure on $G_\mathbb{C}$. See Hall \cite{Hall1994} and references therein.
Give $G$ the bi-invariant metric induced by an $\Ad(G)$-invariant inner product $(\cdot,\cdot)_{\mathfrak{g}}$. Relative to this metric, the tube $G_\epsilon$, of any radius $\epsilon\in (0,\infty]$, is just
\begin{align*}
G_\epsilon = \big\{ \exp(iY) \, \big| \, |Y|_\mathfrak{g}<\epsilon \big\}G.
\end{align*}
The universal property of $G_\mathbb{C}$ 
ensures that every irreducible unitary representation $\xi$ of $G$ extends holomorphically to $G_\mathbb{C}$. Write $ [\xi] \in \widehat{G} $ for its unitary equivalence class, and $d_\xi$ for the dimension of its representation space. Write $\Delta = \Delta_G$ for the Laplacian. The Fourier transform on $G$ is $\mathcal{F}= \mathcal{F}_G$. See also Folland \cite[Chapter 5]{FollandHarmonic} for details.
Recall that components of $\xi$ are eigenfunctions for $-\Delta_G$ for a single eigenvalue $\lambda_\xi$, and $\widehat{G}$ forms a basis for $L^2(G)$ via Peter-Weyl's theorem \cite[Theorem 5.12]{FollandHarmonic}. \\

Combining this, we get our prototypical example, the algebra generated by $\Delta_G$. The Laplacian extends to an operator $(\Delta_G)_\mathbb{C}$ on holomorphic functions on $G_\mathbb{C}$ or $G_\epsilon$.
If $u\in \mathcal{O}(G_\epsilon)$, and $z=\exp(iY)x$ with $x\in G$, we can write
\begin{align*}
(\Delta_G)_\mathbb{C}u(z)
&=
\sum_{[\xi]\in \widehat{G}} d_\xi \Tr (-\lambda_\xi \xi(z)  \mathcal{F}u(\xi)) \\
&=
\sum_{[\xi]\in \widehat{G}} d_\xi \Tr \Big(-\lambda_\xi \xi(x)  \int_G u(\exp(iY)y) \xi(y)^* \, dy \Big),
\end{align*}
which converges absolutely uniformly for $Y$ in compact subsets of $\{Y\in \mathfrak{g} \, | \, |Y|_\mathfrak{g} < \epsilon\}$. Therefore $(\Delta_G)_\mathbb{C}u$ is well-defined, holomorphic on $G_\epsilon$, and $(\Delta_G)_\mathbb{C}$ is bi-$G$-invariant.
If $f\in C^0(\mathbb{C})$ has no more than polynomial growth, we can replace $-\lambda_\xi$ with $f(-\lambda_\xi)$.
For example, if $s\in \mathbb{R}$, we form 
\begin{align*}
(I-(\Delta_G)_\mathbb{C})^\frac{s}{2}u(z) = 
\sum_{[\xi]\in \widehat{G}} d_\xi \Tr (\langle \xi \rangle^s \xi(z)  \mathcal{F}u(\xi)),
\end{align*}
where $\langle \xi \rangle = (1+\lambda_\xi^2)^\frac{1}{2}$,
and this leads to natural "holomorphic Sobolev spaces" on $G_\epsilon$:
\begin{align*}
\{ u\in \mathcal{O}(G_\epsilon) \, | \, 
(I-(\Delta_G)_\mathbb{C})^{\frac{s}{2}}u \in HL^2(G_\epsilon) \}.
\end{align*}
Unfortunately, these ideas do not carry over directly to the Grauert tubes $M_\epsilon$ of $M$. 
How should we define $(I-\Delta_\mathbb{C})^{\frac{s}{2}}$, analogous to $(I-(\Delta_G)_\mathbb{C})^{\frac{s}{2}}$, on $M_\epsilon$ when $s\not\in 2\mathbb{N}_0$?

\newpage
We will overcome this problem by using the Poisson transform due to Stenzel \cite{Stenzel2015}. This transform relates the holomorphic functions on $M_\epsilon$ to their restrictions on $M$. It is used by Stenzel to get an analogue of the Segal-Bargmann transform \cite{Hall1994} for $M_\epsilon$. Arising from a Fourier integral operator, it has well-understood mapping properties, and an entire tube is not required to get a Hilbert space structure on its image \cite{Stenzel2015}.
The Poisson transform relies on an instance of a theorem of Boutet de Monvel \cite{BoutetDeMonvel1978}. This theorem has also been used in several contexts by Zelditch, see e.g. \cite{Zelditch2012, Zelditch2007,Zelditch2014Ergodicity}, but it has not been proved in the generality originally stated by Boutet de Monvel. So we refer to the general statement as a conjecture, although in \cite{BoutetDeMonvel1978} it is a theorem. 
It should be said that a very rough sketch is provided in \cite{BoutetDeMonvel1978}.\\

On $M$ the conjecture implies a supply of operators with the $\epsilon$-extension property. If it holds for a $P\in \PsiPHG^d(M)$, elements in the $\Psi(M)$-commutant of $P$ are extendible,
and their extensions mirror the Sobolev mapping property \cite[pp. 69, Theorem 8.2]{Shubin2001}.
In particular, real-analytic elliptic, and formally normal differential operators have it,
provided that the principal symbol also satisfies a certain strict convexity condition.
This is contained in Theorems \ref{thm:PAlgebraExtendible} and \ref{thm:parametrixresult}, and the two Corollaries \ref{cor:laplacianparametrixresult} and \ref{cor:LaplacianAlgebraExtendible}.
But without a full proof, it is only certain for $P= -\Delta_g$. \\

On the compact Lie group $G$, we set up an even larger algebra in Definition~\ref{def:Sepsilon}. 
This is accomplished by exploiting the fact that $\Psi(G)$ admits a bijection $p\mapsto \Op(p)$, where $p\in S^d(G\times \widehat{G})=S^d$ are symbols due to Ruzhansky, Turunen and Wirth \cite{RuzhanskyTurunenWirth2014}.
These operators are the quantizations of holomorphic symbols in a subalgebra of  $S^d$, and  Theorem~\ref{thm:realansymbelliptic} characterizes ellipticity in it, resembling \cite[pp. 10, Theorem 4.1]{RuzhanskyTurunenWirth2014}. Also here, the Sobolev mapping property is mirrored, in Theorem~\ref{thm:LieGroupSobolevMirror}.

\medskip

\section{Notation}
We write $B(H)$ for the bounded operators on a normed space $H$, and $B(H_1, H_2)$ for operators bounded from one normed space $H_1$ into another $H_2$.
$1_S$ denotes the indicator of a (sub)set $S$. In multiple contexts, $I$ is the identity map. 
Fix a smooth positive 1-density $\omega_0$ on $M$, and form the $L^2$-space on $M$ relative to it.
Then $H^s(M)$ is the Sobolev space of order $s\in \mathbb{R}$ on $M$, normed by
\begin{align*}
||u||_{H^s(M)} = ||(I-\Delta_g)^\frac{s}{2}u ||_{L^2(M)}
\quad
\textnormal{for any}
\quad
u\in H^s(M),
\end{align*}
where $-\Delta_g$ is the positive Laplacian, $(I-\Delta_g)^\frac{s}{2}$ is defined as in \cite[Theorem 10.1]{Shubin2001}.
On a K\"ahler manifold-with-boundary $N$, boundary $\partial N$ is automatically orientable, 
and it is equipped with the volume form induced by the K\"ahler volume form on $N$.
The Sobolev spaces $H^s(N^\circ)$ and $H^s(\partial N)$ are then defined relative to these choices, and the space $\mathcal{O}^s(\partial N)$ is the closure in $H^s(\partial N)$ of 
\begin{align*}
\{ u|_{\partial N} \, | \, u\in C^\infty(\overline{N}) \} \cap \mathcal{O}(N^\circ)
\end{align*}
The distribution space is denoted $\mathcal{D}'(M)$, or $\mathcal{D}'(N^\circ, \Omega^{p,q})$ for distribution $(p,q)$-forms. Finally,
$\partial$ and $\overline{\partial}$ are the Dolbeault operators on $N^\circ$.\\ 

Beyond notation already introduced for $G$, $dx$ is the normalized Haar measure, and $\mathfrak{gl}(k,\mathbb{K})$ is the Lie algebra of 
$k \times k$ matrices with entries contained in a field $\mathbb{K}$. The associated general linear group is
 $\mathrm{GL}(k,\mathbb{K})$, with entries in $\mathbb{K}$. Here $\mathbb{K}$ is $\mathbb{C}$ or $\mathbb{R}$.
Also, $\mathrm{U}(m)$ are the $m\times m$ unitary matrices, and $\mathrm{O}(m)$ the real orthogonal matrices. On matrices, $||\cdot||$
is the norm induced by the Euclidean 2-norm, and $\Tr$ is the trace.

\newpage
\section{Preliminaries} 
Before proceeding, we need some definitions and theorems.
Mainly, theory by Ruzhansky, Turunen and Wirth \cite{RuzhanskyTurunenWirth2014}, Wirth and Ruzhansky \cite{RuzhanskyWirth2014}, and some minor lemmas in Winterrose \cite{Winterrose2021} derived from \cite{RuzhanskyWirth2014}. Let $d\in \mathbb{R}$.\\

Analogous to the Schwartz functions on $\mathbb{R}^n$, there exists a sequence space $\mathcal{S}(\widehat{G})$: 
Define spaces $\mathcal{S}^d(\widehat{G})$ of functions $a$ on the unitary representations $\xi$, such that
\medskip
\begin{enumerate}

\item $a(\xi) \in \mathfrak{gl}(d_\xi, \mathbb{C})$.

\medskip

\item $a(U\xi U^*) = Ua(\xi)U^*$ given any $U \in \mathrm{U}(d_\xi)$.

\medskip

\item $\sum_{[\xi]\in \widehat{G}} d_\xi \langle \xi \rangle^{2d} \Tr(a(\xi)^* a(\xi))  < \infty$.

\end{enumerate}
\medskip

Note that it descends to $a : \widehat{G} \to \mathfrak{gl}(d_\xi, \mathbb{C})/{\sim}$, where $\sim$ is unitary equivalence.
These spaces are equipped with the inner products
\begin{align*}
(a,b)_{\mathcal{S}^d(\widehat{G})} = \sum_{[\xi] \in \widehat{G}} d_\xi \langle \xi \rangle^{2d} \Tr( b(\xi)^* a(\xi))
\quad
\textnormal{for any}
\quad
a,b \in \mathcal{S}^d(\widehat{G}),
\end{align*}
and $\mathcal{S}(\widehat{G}) = \cap_{k\in \mathbb{Z}}\mathcal{S}^k(\widehat{G})$ has the Frechet topology induced by the collection of norms. Then $\mathcal{S}(\widehat{G})$ is sequentially dense in every $\mathcal{S}^{k}(\widehat{G})$ by truncating the sequences in $\mathcal{S}^{k}(\widehat{G})$.
It is the right space; $\mathcal{F} : C^\infty(G) \to \mathcal{S}(\widehat{G})$ is a well-defined linear homeomorphism.
Each space is given the pairing
\begin{align*}
\mathcal{S}^{-k}(\widehat{G}) \times \mathcal{S}^k(\widehat{G}) \to \mathbb{C} : (a,b) \mapsto \langle a, b \rangle = \sum_{[\xi] \in \widehat{G}} d_\xi \Tr( a(\xi)b(\xi)),
\end{align*}
and then Riesz' theorem and sequential denseness of $\mathcal{S}(G)$ implies that
\begin{align*}
\mathcal{S}'(\widehat{G})\cong \cup_{k\in \mathbb{Z}}\mathcal{S}^k(\widehat{G}).
\end{align*}
This is understood in the above pairing, which extends to $\mathcal{S}'(\widehat{G})\times \mathcal{S}(\widehat{G})$. 

\begin{definition}
Overloading notation, $\mathcal{F}$ and $\mathcal{F}^{-1}$ extend by duality:
\medskip
\begin{enumerate}

\item Extend $\mathcal{F}$ to $\mathcal{F} : \mathcal{D}'(G) \to \mathcal{S}'(\widehat{G}) : f \mapsto \mathcal{F}f$ by
\begin{align*}
\langle \mathcal{F} f(\xi), a(\xi) \rangle  
=
\langle f(x), \mathcal{F}^{-1}a(x^{-1}) \rangle   
\quad
\textnormal{for all}
\quad
a \in \mathcal{S}(\widehat{G}).
\end{align*}

\medskip

\item Extend $\mathcal{F}^{-1}$ to  $\mathcal{F}^{-1} : \mathcal{S}'(\widehat{G}) \to \mathcal{D}'(G) : a \mapsto \mathcal{F}^{-1} a$ by
\begin{align*}
\langle \mathcal{F}^{-1} a(x), f(x^{-1}) \rangle  
=
\langle a(\xi),  \mathcal{F}f(\xi) \rangle  
\quad
\textnormal{for all}
\quad
f \in C^\infty(G).
\end{align*}

\end{enumerate}
\medskip
\end{definition}

Continuous operators $P : C^\infty(G) \to C^\infty(G)$ decompose by Peter-Weyl's theorem. The Peter-Weyl expansion of $u\in C^\infty(G)$ converges in $C^\infty(G)$, so we have
\begin{align*}
Pu(x)
&=
\sum_{[\xi]\in \widehat{G}} d_\xi \Tr(\xi(x) p(x,\xi) \mathcal{F}u(\xi))
\quad
\textnormal{for all}
\quad x\in G,
\end{align*}
where $p(x,\xi) = \xi(x)^*(P\xi)(x)$ is the "matrix-symbol" of $P$. This defines $P$ uniquely.
It has the general property that
\begin{align*}
p \in C^\infty (G; \mathcal{S}'(\widehat{G})),
\end{align*}
and we will write $P=\Op(p)$. The map $p \mapsto \Op(p)$ is the operator quantization map.

\newpage
In order to describe the symbol space $S^d$, we need to set up some more structure.
Take a faithful representation, $\rho : G \to \textrm{U}(m)$ for some $m\in \mathbb{N}$, onto a Lie subgroup.
Using real matrices to represent $\mathbb{C}$, we have
\begin{align*}
\textrm{U}(m) \cong \textrm{O}(2m) \cap\textrm{GL}(m,\mathbb{C}) \subset \textrm{GL}(2m,\mathbb{R}),
\end{align*}
and so we may assume $\rho : G \to \textrm{O}(2m) $.
Choose a complement to $d\rho(\mathfrak{g})$ in $\mathfrak{gl}(2m,\mathbb{R})$. By the inverse function theorem, there is an open $V\subset \mathrm{GL}(2m,\mathbb{R})$ with 
\begin{align*}
\rho(G) \subset V,
\end{align*}
and a neighbourhood $ U $ of $0$ in the complement, giving a diffeomorphism
\begin{align*}
G \times U \to V : (g, Y) \mapsto \rho(g)\exp(Y),
\end{align*}
which in turn gives a smooth "projection" map
\begin{align*}
\varrho : V \to G : \rho(g)\exp(Y) \mapsto g.
\end{align*}
It absorbs $G$ from the left, that is,
\begin{align*}
x\varrho(y)=\varrho(\rho(x)y)
\quad
\textnormal{for any}
\quad
(x,y) \in G \times V.
\end{align*} 
Fix a cutoff $\chi\in C^\infty_0(V)$ equal to $1$ on a $\rho(G)$-invariant neighbourhood of $\rho(G)$ in $V$. Note that $\chi(\rho(x)y)=\chi(\rho(x))=1$ for $x\in G$ and $y$ in this neighbourhood. \\

This choice of embedding of $G$ into a Euclidean space is fixed once and for all. 
This is the structure implicit in Ruzhansky, Turunen and Wirth \cite{RuzhanskyTurunenWirth2014}. 

\begin{definition}
Associated to these choices, we define two operators $\delta^\alpha_x$ and $\delta^\alpha_\xi$.
Given any $\alpha\in \mathbb{N}^{m \times m}_0$ with $m$ as in the embedding $\rho : G \to \textrm{O}(2m)$, put
\begin{align*}
d_\alpha(x)=(\rho(x^{-1})-I)^\alpha
\quad
\textnormal{for all}
\quad
x \in G,
\end{align*}
where $\alpha$ is viewed as a multi-index, the power is a product over terms indexed by $\alpha$. 

\medskip
\begin{enumerate}

\item Define $\delta^\alpha_x : C^\infty(G) \to C^\infty(G) : u \mapsto \delta^\alpha_x u $ by
\begin{align*}
\delta^\alpha_x u(x) &= \big[ \partial^\alpha_y [ \chi(u\circ \varrho)(\rho(x)y)] \big] \big|_{y=I}
\quad
\textnormal{for all}
\quad
x\in G.
\end{align*}

\medskip

\item Define $\delta^\alpha_\xi : \mathcal{S}'(\widehat{G}) \to \mathcal{S}'(\widehat{G}) : a \mapsto \delta^\alpha_\xi a $ by
\begin{align*}
\delta^\alpha_\xi a(\xi) = \mathcal{F} d_\alpha \mathcal{F}^{-1} a
\quad
\textnormal{for all}
\quad
[\xi] \in \widehat{G}.
\end{align*}

\end{enumerate}
\medskip
\end{definition}

The "difference operators" $\delta_\xi^\alpha$ satisfy a Leibniz-type rule on the $\mathcal{S}'(\widehat{G})$ sequences. 
Write $e_{ij}$ for the $m \times m$ matrix with $1$ at the $i,j$ entry and zero else.

\begin{lemma}[Ruzhansky, Turunen and Wirth \cite{RuzhanskyTurunenWirth2014}]  \label{lmm:DifferenceOperatorLeibnizRule}
\begin{align*}
\delta^{e_{ij}}_\xi(ab)
=
\delta^{e_{ij}}_\xi (a)b
+
a\delta^{e_{ij}}_\xi (b)
+
\sum_{k=1}^m \delta^{e_{ik}}_\xi (a) \delta^{e_{kj}}_\xi (b)
\quad
\textnormal{for any}
\quad
a,b \in \mathcal{S}'(\widehat{G}).
\end{align*}
\end{lemma}

\newpage
Let $d \in \mathbb{R}$. Let us write $X^\beta = X_1^{\beta_1} \circ \cdots \circ X_n^{\beta_n}$ for any $\beta\in \mathbb{N}_0^n$ if $\{X_j\}_{j=1}^n \subset\mathfrak{g}$.

\begin{definition}
Write $p\in S^d(G \times \widehat{G})= S^d$ for $d\in \mathbb{R}$, if the following holds:  
\medskip
\begin{enumerate}

\item $p \in C^\infty (G; \mathcal{S}'(\widehat{G}) ) $.

\medskip

\item Given any ordered basis $X=(X_1, \cdots , X_n)$ of $\mathfrak{g}$, we have 
\begin{align*}
\sup_{(x,[\xi])\in G \times \widehat{G}}\langle \xi \rangle^{|\alpha|-d} || \delta^\alpha_\xi X^\beta_x p(x,\xi) || < \infty
\quad
\textnormal{for any}
\quad
\alpha\in \mathbb{N}_0^{m \times m}
\quad
\textnormal{and}
\quad
\beta\in \mathbb{N}_0^n.
\end{align*}
\end{enumerate}
\medskip
The space $S^d$ has the Frechet topology induced by the above collection of semi-norms. Given a sequence $\{p_j\}_{j=0}^\infty$ with $p_j \in S^{d_j}$ and $d_j \to -\infty$ as $j\to \infty$, we write
\begin{align*}
p \sim \sum_{j=0}^\infty p_j
\quad
\textnormal{if}
\quad
p- \sum^{k-1}_{j=0} p_j \in S^{\max_{j\geq k} d_{j}}
\quad
\textnormal{for each}
\quad
k\in \mathbb{N}.
\end{align*}
These matrix-symbols are called the H\"ormander symbols.
\end{definition}


\begin{theorem}[Ruzhansky, Turunen and Wirth \cite{RuzhanskyTurunenWirth2014}]
\begin{align*}
\Psi^d(G)=\Op\, S^d.
\end{align*}
\end{theorem}


\begin{definition}
Let $d_1,d_2 \in \mathbb{R}$. Given any $p\in S^{d_1}$ and $q \in S^{d_2}$, define:
\medskip
\begin{enumerate}

\item Pointwise in $(x,[\eta])\in  G \times \widehat{G}$ a symbol $p\odot q$ by
\begin{align*}
(p \odot q)(x,\eta)
=
\sum_{[\xi] \in \widehat{G}}
d_\xi \int_G  \Tr\Big(\xi(y^{-1} x) p(x,\xi)\Big) \eta(x^{-1}y) q(y,\eta)   \, dy.
\end{align*}
\medskip

\item Pointwise in $(x,[\eta])\in  G \times \widehat{G}$ a symbol $p^\dagger$ by
\begin{align*}
p^\dagger(x,\eta)
=
\sum_{[\xi] \in \widehat{G}}
d_\xi \int_G  \Tr\Big(\xi(y^{-1} x) p(y,\xi)^* \Big) \eta(x^{-1}y) \, dy.
\end{align*}
\end{enumerate}
\medskip
\end{definition}

Both of the two sums in the above definition are uniformly absolutely convergent. To see this, integrate by parts entry-wise with $I-\Delta_y$ to summon powers of $\langle \xi\rangle^{-1}$. 
Then the sums are dominated by convergent sums scaling with $\langle \eta \rangle$, independent of $x$, because $p$ and its derivatives in $x$ grow polynomially in $\langle \xi\rangle$ at a fixed rate. 

\begin{theorem}[See e.g. Winterrose \cite{Winterrose2021}]
Let $d_1,d_2 \in \mathbb{R}$. The following holds:
\medskip
\begin{enumerate}

\item The map $S^{d_1} \times S^{d_2} \to S^{d_1+d_2} : (p,q) \mapsto p \odot q$ is well-defined, and
\begin{align*}
\Op(p)\Op(q) = \Op(p \odot q)
\quad
\textnormal{with}
\quad
p \odot q \sim \sum_{N=0}^\infty \sum_{|\alpha|=N} \frac{1}{\alpha!} \delta_\xi^\alpha (p) \delta_x^\alpha (q).
\end{align*}

\medskip

\item The map $S^{d_1} \to S^{d_1} : p \mapsto p^\dagger$ is well-defined, and
\begin{align*}
\Op(p)^* = \Op(p^\dagger)
\quad
\textnormal{with}
\quad
p^\dagger \sim \sum_{N=0}^\infty \sum_{|\alpha|=N} \frac{1}{\alpha!} \delta_\xi^\alpha \delta_x^\alpha (p^*).
\end{align*}

\end{enumerate}
\medskip
Both maps are continuous with respect to the Frechet topologies on the symbol spaces.
\end{theorem}

\newpage
In our applications to compact Lie groups, we need the following two lemmas. Both are proved in Winterrose \cite{Winterrose2021}, but are just easy adaptations of proofs in \cite{RuzhanskyWirth2014}, where a parameter $Y\in K$ running in a topological space $K$ has just been inserted. It represents $Y\in \mathfrak{g}$ in the Cartan decomposition $\exp(iY)x=z\in G_\mathbb{C}$ later.

\begin{lemma}[Winterrose \cite{Winterrose2021}] \label{lmm:LieGroupCalderonVaillancourt}
Let $K$ be compact, and $\{ p_Y \}_{Y\in K} \subset S^d$ bounded.
If all $x$-derivatives of $p_Y(x,\xi)$ are continuous in $(x,Y)$, then for $s\in \mathbb{R}$ we have
\begin{align*}
\sup_{Y \in K} || \Op(p_Y) ||_{B(H^{s}(G), H^{s-d}(G))} < \infty,
\end{align*} 
which is to say, $\Op(p_Y)$ is bounded uniformly in $Y$ from $H^{s}(G)$ to $H^{s-d}(G)$.
\end{lemma}

\begin{lemma}[Winterrose \cite{Winterrose2021}] \label{lmm:invmatsymb}
Let $K$ be compact, and $\{ p_Y \}_{Y\in K} \subset S^d$ bounded.
Assume that each $p_Y(x,\xi)$ is point-wise invertible as a matrix belonging to $\mathfrak{gl}(d_\xi, \mathbb{C})$. 
If all $x$-derivatives of $p_Y(x,\xi)$ are continuous in $(x,Y)$, then
\begin{align*}
\sup_{Y\in K}\sup_{(x,[\xi])\in G \times \widehat{G}}\langle \xi \rangle^{d} ||p_Y (x,\xi)^{-1} ||  < \infty,
\end{align*}
and $\{ p_Y^{-1} \}_{Y\in K} \subset S^{-d}$ is bounded, with the same continuity property.
\end{lemma}

\medskip
Recall the calculus of poly-homogeneous (PHG) operators $\PsiPHG^{d}(M)$ in  $\Psi^{d}(M)$. These have a unique homogeneous principal symbol, the classical principal symbol, and elliptic $A\in \PsiPHG^{d}(M)$ with $d>0$ have either discrete spectrum $\sigma(A)$, or $\sigma(A)= \mathbb{C}$.
Choose such an $A\in \PsiPHG^d(M)$, but with a discrete spectrum in the right half-plane.
It is parameter elliptic w.r.t. a sector $\Lambda$ if the principal symbol does not valuate in $\Lambda$.
Suppose $A$ is parameter-elliptic w.r.t. a $[\pi-\theta_0,\pi+\theta_0]$-sector for some $\theta_0 \in (0,\pi)$.
Given $k\in \mathbb{N}_0$ and $z\in \mathbb{C}$ with $\Real(z)-k<0$, we can unambiguously define
\begin{align*}
A^z = \frac{1}{2\pi i} \int_{\Gamma_R} \lambda^{z-k} (\lambda I - A)^{-1} A^k \, d\lambda,
\end{align*}
where $\Gamma_R$ is the counter-clockwise relative to $ \sigma(A)\setminus \{0\}$ boundary of
\begin{align*}
\Lambda_R = \cup_{\theta \in [-\pi+\theta_0,\pi-\theta_0]} e^{i\theta}[R,\infty) \supset \sigma(A)\setminus \{0\},
\end{align*}
and the integral converges in $B(H^N(M))$ for all $N\in \mathbb{N}$ by 
\cite[pp. 86, Theorem 9.3]{Shubin2001}. By the Sobolev embedding theorem, it converges weakly in the topology of $C^\infty(M)$.
Suppose now $A$ is parameter elliptic w.r.t. a $[\theta_0', 2\pi - \theta_0']$-sector for some $\theta_0' \in (0,\frac{1}{2}\pi)$. In the same way, for $\Real(z)>0$, we can form
\begin{align*}
e^{- t A^z}
=
\frac{1}{2\pi i} \int_{\Gamma_R'} e^{- t \lambda^z} (\lambda I - A)^{-1} \, d\lambda
+
\frac{1}{2\pi i} \int_{R \mathbb{S}^1} (\lambda I - A)^{-1} \, d\lambda,
\end{align*}
where $R \mathbb{S}^1$ is oriented counter-clockwise, $\Gamma_R' = \partial \Lambda_R' $ with $R>0$ such that
\begin{align*}
\Lambda_R' = \cup_{\theta' \in [-\theta_0' , \theta_0']} e^{i\theta'}[R,\infty) \supset \sigma(A)\setminus \{0\},
\end{align*}
and $\Gamma_R'$ is oriented counter-clockwise relative to $ \sigma(A)\setminus \{0\}$ inside the right half-plane. In particular, if $A$ is formally self-adjoint and elliptic, its spectrum must be discrete, and is semi-bounded from below if its classical principal symbol is positive on $T^*M\setminus 0$. This is \cite[pp. 71, Theorem 8.3-8.4]{Shubin2001} and  \cite[pp. 86, Corollary 9.3]{Shubin2001}.

\newpage
\section{The Boutet de Monvel Conjecture} (A theorem in a prototypical case.)
A theorem was announced by Boutet de Monvel in 
an old 1978 conference paper \cite{BoutetDeMonvel1978}, but only the case $\sqrt{-\Delta_g}$ in \cite{BoutetDeMonvel1978} was proved by Stenzel in \cite{Stenzel2014} and Zelditch in \cite{Zelditch2012}.
Several possible approaches to proving it generally are discussed by Zelditch in \cite{Zelditch2012}.
To our knowledge, the general statement is not yet fully proven.\\ 

Let $d>0$, and let $P\in \PsiPHG^d(M)$ be analytic with classical principal symbol $p$. 
Suppose that $P$ is formally self-adjoint, elliptic, $p|_{T^*M\setminus 0} > 0$, and that
\begin{align*}
\{ \xi \in T^*_x M \, | \, p(x,\xi) \leq 1 \}
\quad \textnormal{is strictly convex for all}
\quad
x\in M.
\end{align*}
In that case, $P$ is semi-lower bounded, and we assume it has no negative eigenvalues.
A real-analytic complete Hamiltonian flow
$\varphi_t : T^* M \setminus 0 \to T^*M \setminus 0 $ arises from $p^\frac{1}{d}$. Given any $t\in \mathbb{R}$ and $(x,\xi)\in T^*M\setminus 0$, it satisfies
\begin{align*}
\varphi_{t}(x,s\xi) = s\cdot \varphi_{t}(x,\xi)
\quad
\textnormal{for all}
\quad
s>0,
\end{align*}
and $\varphi_t(x,\xi)$ extends holomorphically in $t$ to $\{ t\in \mathbb{C} \, | \, |t|<\epsilon \}$ for $\epsilon>0$ small enough.
The extension is then a map
\begin{align*}
\{ t\in \mathbb{C} \, | \, |t|<\epsilon \} \to T^*_{1,0}M_\mathbb{C} \setminus 0 : t \mapsto \varphi_t(x,\xi),
\end{align*}
where $T^*_{1,0}M_\mathbb{C}$ is the holomorphic cotangent bundle of $M_\mathbb{C}$.\\

\begin{conjecture}[Boutet de Monvel \cite{BoutetDeMonvel1978}. The statement here is a bit different] \label{conj:BoutetDeMonvelConjecture}
Let $P_\epsilon$ denote the Schwartz kernel of $e^{-\epsilon P^\frac{1}{d}}$, and put
\begin{align*}
\Phi_x(\xi) = (\pi_x \varphi_{i p(x,\xi)^\frac{1}{d}}) (x, \xi)
\quad
\textnormal{for any}
\quad
(x,\xi)\in B^*_\epsilon M \setminus 0,
\end{align*}
where $\pi_x$ is the bundle projection, and
\begin{align*}
B^*_\epsilon M = \{ (x,\xi) \in T^* M \, | \, p(x,\xi)^\frac{1}{d} < \epsilon  \}.
\end{align*}
There is a maximal $\epsilon_0>0$ such that for any $\epsilon\in (0,\epsilon_0)$ the following holds: 
\medskip
\begin{enumerate}

\item The extended flows combine into a real-analytic diffeomorphism
\begin{align*}
\Phi : B_\epsilon^*M\setminus 0 \to M_\epsilon^\Phi \setminus M \subset M_\mathbb{C} : (x,\xi) \mapsto \Phi_x(\xi).
\end{align*}

\medskip

\item The set $M_\epsilon^\Phi$ is open in $M_\mathbb{C}$ with (orientable) $C^\omega$-boundary $\partial M_\epsilon^\Phi$.

\medskip

\item The set $M_\epsilon^\Phi$ is strictly pseudo-convex.

\end{enumerate}
\medskip
(Point 3. implies that $M_\epsilon^\Phi$ admits a K\"ahler metric $h$, and a global potential for it.)
There is a maximal $\epsilon'_0 \in (0, \epsilon_0]$ such that for any $\epsilon \in (0,\epsilon_0')$ the following holds:
\medskip
\begin{enumerate}
\setcounter{enumi}{3}

\item The map $x \mapsto P_\epsilon(x,y)$ extends holomorphically to $M_\epsilon^\Phi$ for each fixed $y\in M$.

\medskip

\item The kernel $P_\epsilon |_{\partial M_\epsilon^\Phi \times M}$ induces a complex-phase PHG FIO $S_\epsilon$ of order $-\frac{n-1}{4}$.

\medskip

\item $S_\epsilon$ defines a homeomorphism $S_\epsilon : H^s(M) \to \mathcal{O}^{s+\frac{n-1}{4}}(\partial M_\epsilon^\Phi)$ for any $s\in \mathbb{R}$. 

\end{enumerate} 
\medskip
\end{conjecture}

The first three points  simply generalize the basic properties of the Grauert tubes.
In the special case that $P=-\Delta_g$, $M_\epsilon^\Phi$ is exactly the Grauert tube $M_\epsilon$.









\newpage
The abbreviation "PHG FIO" is for poly-homogeneous Fourier integral operator. In this article, we only remind the reader of the classical article by H\"ormander \cite{FIO1}, and, in the situation of globally defined complex phase FIO, Melin and Sj\"ostrand \cite{ComplexFIO}.
The theory in the complex case runs parallel to the real case, but is less widely known. 
In the proven instance of \ref{conj:BoutetDeMonvelConjecture}, there are precise details on the structure of $S_\epsilon$:\\

\begin{theorem}[Stenzel \cite{Stenzel2014, Stenzel2015} and Zelditch \cite{Zelditch2012}. See Boutet de Monvel \cite{GuilleminBoutetDeMonvel1981}] \label{thm:BoutetDeMonvelTheorem}
Conjecture~\ref{conj:BoutetDeMonvelConjecture} is true for the Laplacian $P=-\Delta_g$ of $(M,g)$ for some $\epsilon_0\geq \epsilon_0'>0$. The positive homogeneous canonical relation of $S_\epsilon$ arises from the graph of 
\begin{align*}
T^*M\setminus 0 \to T^*(\partial M_\epsilon)\setminus 0 : (x,\xi) \mapsto \frac{|\xi|_x}{\epsilon} \alpha_{(x,\xi)}, 
\end{align*}
where the cotangent vector is
\begin{align*}
\alpha_{(x,\xi)} = (\iota^*_{\partial M_\epsilon}\Imag \overline{\partial} \rho)_{\exp_x (i\epsilon |\xi|_x^{-1}\xi^\sharp)},
\end{align*}
and $\rho$ is the unique tube function of $M_\epsilon$, a K\"ahler potential for the structure on $M_\epsilon$.
This function is the unique real-analytic solution $\rho : M_\epsilon \to [0,\infty)$ of  
\begin{align*}
\begin{cases}
(i\partial \overline{\partial} \sqrt{\rho})^{\wedge n} = 0 & \textnormal{in} \quad M_\epsilon \setminus M,  \\
\textnormal{Lev}(\rho)|_M = \frac{1}{2} g, & 
\end{cases}
\end{align*}
where $\textnormal{Lev}(\rho)$ is the Levi form of $\rho$. The second condition ensures isometric inclusion.
It has, and for uniqueness must have, the following properties:
\medskip
\begin{enumerate}

\item $\rho$ is invariant under the unique anti-holomorphic involution on $M_\epsilon$.

\medskip

\item $\rho$ is strictly plurisubharmonic on $M_\epsilon$.

\medskip

\item $\rho|_M=0$ and $d\rho|_M = 0$.

\end{enumerate}
\medskip
Finally, $S^*_\epsilon S_\epsilon \in \PsiPHG^{-\frac{n-1}{2}}(M)$ is elliptic, with classical principal symbol
\begin{align*}
T^*M \to [0,\infty) : (x,\xi) \mapsto |\xi^\sharp|_x^{-\frac{n-1}{2}}.
\end{align*}
\end{theorem}

\medskip
\begin{example}[Torus $\mathbb{T}^n$]
{\rm 
Equip the $n$-torus $\mathbb{T}^n \cong (\mathbb{R}/2\pi\mathbb{Z})^n $ with the flat metric. The Bruhat-Whitney complexification of $\mathbb{T}^n$ is 
\begin{align*}
\mathbb{T}^n_{\mathbb{C}} \cong (\mathbb{C}/2\pi\mathbb{Z})^n,
\end{align*}
and the Grauert tube function of $\mathbb{T}^n_\epsilon$ for any $\epsilon>0$ is
\begin{align*}
\sqrt{\rho} : \mathbb{T}^n_{\epsilon}  \to [0,\infty) : [z] \mapsto   |\Imag (z)|,
\end{align*}
where $\mathbb{T}^n_\epsilon \cong (\mathbb{R}_\epsilon /2\pi \mathbb{Z})^n$ for $\mathbb{R}_\epsilon = \{ z\in \mathbb{C} \, | \, |\Imag(z)<\epsilon \}$, $2\pi \mathbb{Z}$ acting on the real part.
The half-wave propagator acting on $u\in C^\infty(\mathbb{T}^n)$ is then just
\begin{align*}
e^{-\epsilon\sqrt{-\Delta_{\mathbb{T}^n}}}u(x) =
\sum_{k \in \mathbb{Z}} e^{ik\cdot x - \epsilon |k|} \int_{\mathbb{T}^n} e^{-ik\cdot y} u(y) \, dy,
\end{align*}
where the kernel extends to $x\in \mathbb{T}^n_\epsilon$, but is singular for $|x|=\epsilon$.
}
\end{example}


\newpage
\section{Operators with the $\epsilon$-Extension Property}
Let $n = \dim M$ throughout.
The Poisson transform in \cite{ Stenzel2015} is built from the solution operator to
\begin{align*}
\begin{cases}
\Delta_h u = 0 & \textnormal{in} \quad M_\epsilon, \\
u|_{\partial M_\epsilon} = f & \textnormal{on} \quad \partial M_\epsilon,
\end{cases}
\end{align*}
where $\Delta_h$ is the K\"ahler Laplacian associated to the K\"ahler metric $h$ from Theorem~\ref{thm:BoutetDeMonvelTheorem}, and $f$ is, at first, a smooth function on $\partial M_\epsilon$, and we look for $u\in C^2(M_\epsilon) \cap C^0(\overline{M_\epsilon})$. The unique solvability of this Dirichlet problem for the Laplacian $\Delta_h$ is well-known, but we need details on extensions of the solution operator. \\

Assume therefore that $(N, h)$ is a smooth Riemannian manifold-with-boundary. In our case, $N$ is the closure of an open subset $N^\circ$ in a compact smooth manifold $\widetilde{N}$, where $N=N^\circ \sqcup \partial N$, the topological boundary is a co-dimension one hypersurface. 
(In further applications, either $N=M_\epsilon$ or $N= M^\Phi_\epsilon$, carrying the K\"ahler metric $h$.)
This "enveloping manifold" $\widetilde{N}$ carries a metric coinciding with $h$ on $N$.

\begin{definition}
Given any $s\in \mathbb{R}$, we define the $H^s$-extendible distributions
\begin{align*}
\bar{H}^s(N^\circ)
=
\{ u|_{N^\circ} \in \mathcal{D}'(N^\circ) \, | \, u\in H^s(\widetilde{N}) \},
\end{align*}
and equip this space with the norm
\begin{align*}
||v||_{\bar{H}^s(N^\circ)}
=
\inf_{u|_{N^\circ} = v} ||u||_{H^s(\widetilde{N})}
\quad
\textnormal{if}
\quad
v\in \bar{H}^s(N^\circ). 
\end{align*}
\end{definition}

These spaces of $s\in \mathbb{R}$ extendible distributions are ingredients in the next theorem. It states that the Dirichlet problem for $\Delta_h$ on $N$ is uniquely solvable in these spaces. A clear exposition (of the invertible case) is found in Grubb \cite[pp. 320-326]{GG}.

\begin{theorem}[Outlined in Grubb \cite{GG}. Alternatively, see Boutet de Monvel \cite{BoutetDeMonvel1979}] \label{thm:PoissonSolutionOperator}
Let $s\in \mathbb{R}$. There are bounded, linear, and mutually inverse bijections
\begin{align*}
K_\gamma &: H^{s-\frac{1}{2}}(\partial N)  \to \bar{H}^s (N^\circ) \cap \ker\, (\Delta_h|_{\bar{H}^s(N^\circ)}) , \\
\gamma_0 &: \bar{H}^s(N^\circ) \cap \ker\, (\Delta_h|_{\bar{H}^s(N^\circ)}) \to H^{s-\frac{1}{2}}(\partial N).
\end{align*}
If $f\in H^{s-\frac{1}{2}}(\partial N)$, they together uniquely solve
\begin{align*}
\begin{cases}
\Delta_h (K_\gamma f) = 0 & \textnormal{in} \quad N, \\
\gamma_0 (K_\gamma f) = f & \textnormal{on} \quad\partial N.
\end{cases}
\end{align*}
The operator $\gamma_0$ acts as the restriction operator onto $\partial N$ for smooth functions on $\overline{N}$,
and $K_\gamma$ is a Poisson operator (see \cite{GG}), restricting to a continuous linear map 
\begin{align*}
K_\gamma : C^\infty(\partial N) \to C^\infty(N).
\end{align*}
Relative to fixed smooth positive 1-densities on $N$ and $\partial N$, it has the two properties:
\medskip
\begin{enumerate}

\item $K_\gamma$ admits a unique formal $L^2$-adjoint $K^*_\gamma$.

\medskip

\item  $K^*_\gamma K_\gamma \in \PsiPHG^{-1}(\partial N)$ has classical principal symbol positive on $T^*N \setminus 0$.

\end{enumerate}
\medskip
\end{theorem}

\newpage
In the case that $(N^\circ,h)$ is K\"ahler, solutions produced by $\mathcal{O}^s(\partial N)$ are holomorphic.
Naturally, the boundary $\partial N$ carries the Riemannian structure induced from $N$.

\begin{lemma} \label{lmm:PoissonImage}
In the above setting, if $(N^\circ,h)$ is K\"ahler, we have
\begin{align*}
K_\gamma \mathcal{O}^s(\partial N) \subset \Big\{ \, v\in \bar{H}^{s+\frac{1}{2}}(N^\circ) \cap C^\infty(N^\circ) \, \Big| \, \overline{\partial} (v|_{N^\circ}) = 0 \, \Big\}.
\end{align*}
\end{lemma}
\begin{proof}
Pick any element $u \in \mathcal{O}^s(\partial N)$. We must then show that $\overline{\partial} ( K_\gamma u |_{N^\circ} ) = 0$.
Take a sequence $\{u_j\}_{j=1}^\infty \subset C^\infty(\overline{N})$ with $u_j \in \mathcal{O}(N^\circ)$ such that
\begin{align*}
u_j|_{\partial N} \to u
\quad
\textnormal{in}
\quad
H^s(\partial N)
\quad
\textnormal{as}
\quad
j \to \infty.
\end{align*}
Since holomorphic functions are harmonic for $\Delta_h$, we have
\begin{align*}
u_j|_{N^\circ} \in \ker\, (\Delta_h|_{\bar{H}^s(N^\circ)}).
\end{align*}
Consequently, we can write
\begin{align*}
u_j=K_\gamma(u_j|_{\partial N}) \to K_\gamma u 
\quad
\textnormal{in}
\quad
\bar{H}^{s+\frac{1}{2}}(N^\circ) 
\quad
\textnormal{as}
\quad
j \to \infty,
\end{align*}
and then also
\begin{align*}
0=\overline{\partial} (u_j |_{N^\circ})  \to \overline{\partial} ( K_\gamma u |_{N^\circ} )
\quad
\textnormal{in}
\quad
\mathcal{D}'(N^\circ, \Omega^{0,1})
\quad
\textnormal{as}
\quad
j \to \infty.
\end{align*}
Thus $\overline{\partial} ( K_\gamma u |_{N^\circ} ) = 0$,  and $K_\gamma u \in C^\infty(N^\circ)$ by ellipticity of $\Delta_h$. 
\end{proof}


Choose $\epsilon_0'\leq \epsilon_0$ so that the Boutet de Monvel Theorem~\ref{thm:BoutetDeMonvelTheorem} holds, and $\epsilon\in(0,\epsilon_0')$. Using Theorem~\ref{thm:PoissonSolutionOperator}, we construct the Poisson transform as defined by Stenzel in \cite{Stenzel2015}. The Poisson transform is well-behaved, and can be unitarized by polar decomposition. It is very closely related to the restriction map \cite{Stenzel2015}. \\

Let $h$ be the K\"ahler metric on the closure of a tube $M_\epsilon$, and let $s\in \mathbb{R}$ be arbitrary.

\begin{definition}
If $f\in H^{s-\frac{n+1}{4}}(M)$, $\mathcal{P}_\epsilon f$ is the extension of $e^{-\epsilon \sqrt{-\Delta_g}}f$ to $M_\epsilon$.
\end{definition}

Suppose $K_\epsilon$ is the Poisson operator solving the Dirichlet problem for $\Delta_h$ on $\overline{M_\epsilon}$. Then, by Theorem~\ref{thm:PoissonSolutionOperator} it is a linear homeomorphism
\begin{align*}
K_\epsilon :H^{s-\frac{1}{2}}(\partial M_\epsilon)  \to \bar{H}^s (M_\epsilon) \cap \ker(\Delta_h |_{\bar{H}^s (M_\epsilon)}),
\end{align*}
and so, by Theorem~\ref{thm:BoutetDeMonvelTheorem}, $\mathcal{P}_\epsilon= K_\epsilon S_\epsilon$ realizes a linear homeomorphism 
\begin{align*}
\mathcal{P}_\epsilon : H^{s-\frac{n+1}{4}}(M) \to K_\epsilon \mathcal{O}^{s-\frac{1}{2}}(\partial M_\epsilon).
\end{align*}
In the special case that $s=0$, we have 
\begin{align*}
K_\epsilon \mathcal{O}^{-\frac{1}{2}}(\partial M_\epsilon) = HL^2(M_\epsilon),
\end{align*}
which follows from $\bar{H}^{0}(M_\epsilon) = L^2(M_\epsilon)$, as any $L^2$-function can be extended by zero, and a formulation of $\mathcal{O}^{-\frac{1}{2}}(\partial M_\epsilon)$ as the kernel of the co-boundary operator $\overline{\partial}_b$ on $\partial M_\epsilon$. See the remarks in \cite{Stenzel2014, Stenzel2015, Zelditch2012}.

\begin{theorem}
The composite $\mathcal{P}_\epsilon : H^{s-\frac{n+1}{4}}(M) \to K_\epsilon \mathcal{O}^{s-\frac{1}{2}}(\partial M_\epsilon)$ is well-defined. Furthermore, it is a linear homeomorphism.
\end{theorem}

\newpage
Our approach is to indirectly define spaces as images under $\mathcal{P}_\epsilon$ of any order $s\in \mathbb{R}$. These will coincide with the natural definition on $G_\epsilon$ when $M=G$.
\begin{definition}
\begin{align*}
HH^s(M_\epsilon) = \mathcal{P}_\epsilon H^{s-\frac{n+1}{4}}(M).
\end{align*}
\end{definition}
Using the inverse $\mathcal{P}_\epsilon^{-1}$, the image space is made to inherit the Hilbert space structure. In the case $M=G$ and $P=\Delta_G$, it has the inner product induced by
\begin{align*}
||u||_{HH^s(G_\epsilon)}
=
||\mathcal{P}_\epsilon (I-\Delta_G)^{\frac{s}{2}} \mathcal{P}_\epsilon^{-1} u ||_{L^2(G_\epsilon)} 
\quad
\textnormal{if}
\quad
u\in HH^s(G_\epsilon),
\end{align*}
and $\mathcal{P}_\epsilon : H^{s-\frac{n+1}{4}}(G) \to HH^s(G_\epsilon)$ is automatically a bounded isomorphism this way. This is clear from 
\begin{align*}
|| \mathcal{P}_\epsilon u ||_{HH^s(G_\epsilon)}
&=
|| \mathcal{P}_\epsilon (I-\Delta_G)^{\frac{s}{2}} u ||_{L^2(G_\epsilon)} \\
&\leq ||\mathcal{P}_\epsilon||_{B(H^{-\frac{n+1}{4}}(G), HL^2(G_\epsilon))} || u ||_{H^{s-\frac{n+1}{4}}(G)}.
\end{align*}

\begin{proposition}
\begin{align*}
HH^s(G_\epsilon)
=
\{ u\in \mathcal{O}(G_\epsilon) \, | \, 
(I-(\Delta_G)_\mathbb{C})^{\frac{s}{2}}u \in HL^2(G_\epsilon) \}.
\end{align*}
\end{proposition}
\begin{proof}
Take any $u\in C^\infty(G)$, and for any $z\in G_\epsilon$ write
\begin{align*}
\mathcal{P}_\epsilon (I-\Delta_G)^\frac{s}{2}u(z)
&=
\sum_{[\xi]\in \widehat{G}} d_\xi \Tr ( \mathcal{P}_\epsilon\xi(z) \mathcal{F} (I-\Delta_G)^\frac{s}{2}u ) \\
&=
\sum_{[\xi]\in \widehat{G}} d_\xi \Tr \Big( \langle \xi \rangle^s 
\xi(z) \int_G u(x)  e^{-\epsilon\sqrt{\lambda_\xi} }\xi(x)^* \, dx \Big) \\
&=
\sum_{[\xi]\in \widehat{G}} d_\xi \Tr \Big( \langle \xi \rangle^s 
\xi(z) \int_G (e^{-\epsilon\sqrt{-\Delta}}u)(x) \xi(x)^* \, dx \Big),
\end{align*}
where we use that the Peter-Weyl expansion of $u$ converges in the topology of $C^{\infty}(G)$, and so the first sum converges in $HL^2(G_\epsilon)$, thus uniformly on compact subsets of $G_\epsilon$.
This shows that
\begin{align*}
\mathcal{P}_\epsilon (I-\Delta_G)^\frac{s}{2}u = (I-(\Delta_G)_\mathbb{C})^\frac{s}{2} \mathcal{P}_\epsilon  u,
\end{align*}
which, on the left hand side, must extend to $u\in H^{s-\frac{n+1}{4}}(G)$ just by the continuity. To see equality for $u\in H^{s-\frac{n+1}{4}}(G)$, pick $\{u_k\}_{k=1}^\infty\subset C^\infty(G)$ such that
\begin{align*}
u_k \to u
\quad
\textnormal{in}
\quad
H^{-\frac{n+1}{4}}(G)
\quad
\textnormal{as}
\quad
k\to \infty,
\end{align*}
and note that
\begin{align*}
(I-(\Delta_G)_\mathbb{C})^\frac{s}{2} \mathcal{P}_\epsilon u_k \to (I-(\Delta_G)_\mathbb{C})^\frac{s}{2} \mathcal{P}_\epsilon u
\quad \textnormal{uniformly on} \quad
K\subset\subset G_\epsilon
\quad \textnormal{as}
\quad k\to \infty, 
\end{align*}
while the same type of convergence is implied by $HL^2$-convergence on the left side.
Therefore, the two limits agree as holomorphic functions on $G_\epsilon$ when $u\in H^{s-\frac{n+1}{4}}(G)$, and
the proposition follows.
\end{proof}

\newpage
\subsection{On a Compact Riemannian Manifold $M$}
Let $M$ carry the metric $g$.
The key is the following simple lemma:

\begin{lemma} \label{lmm:PowersExponentsCommutator}
Let $d\in \mathbb{N}$ and let $P\in \PsiPHG^d(M)$ be elliptic and formally self-adjoint.
Assume that $P$ has classical principal symbol $p$ positive on $T^* M \setminus 0$, and 
\begin{align*}
\sigma(P) \subset [0,\infty).
\end{align*}
Then, if $A : C^\infty(M) \to C^\infty(M)$ is continuous with $[A,P]=0$, the following holds:
\medskip
\begin{enumerate}

\item $[A, P^\frac{1}{d}] = 0$.

\medskip

\item $[A, e^{-\epsilon P^\frac{1}{d}}] = 0$ for any $\epsilon>0$.

\end{enumerate}
\medskip
\end{lemma}


\begin{proof}
By hypothesis, $P$ is parameter-elliptic w.r.t. any closed sector in $\mathbb{C}\setminus(0,\infty)$.
Thus the two operators are defined, given appropriate contours in the complex plane. 
That is, we can choose $R>0$ so that if $u\in C^\infty(M)$, we have
\begin{align*}
P^\frac{1}{d}u  =
\frac{1}{2\pi i}\int_{\Gamma_R} \lambda^\frac{1-d}{d} ( \lambda I - P )^{-1}Pu \, d\lambda,
\end{align*}
and
\begin{align*}
e^{- t P^\frac{1}{d}}u
=
\frac{1}{2\pi i} \int_{\Gamma_R'} e^{- t \lambda^\frac{1}{d}} (\lambda I - P)^{-1}u \, d\lambda
+
\frac{1}{2\pi i} \int_{R \mathbb{S}^1} (\lambda I - P)^{-1}u \, d\lambda,
\end{align*}
where $\lambda \mapsto \lambda^\frac{1-d}{d} $ is defined by the principal logarithm with branch cut along $(-\infty,0]$,
and $\Gamma_R$ and $\Gamma_R'$ are keyhole contours, $R\mathbb{S}^1$ encircles no eigenvalues except possibly $0$.
The commutators are easily calculated to be
\begin{align*}
[A, P^\frac{1}{d}]u
&=
A\Big(\frac{1}{2\pi i}\int_{\Gamma_R} \lambda^\frac{1-d}{d} ( \lambda I - P )^{-1}Pu \, d\lambda \Big)
- P^\frac{1}{d}Au \\
&=
\frac{1}{2\pi i}\int_{\Gamma_R} \lambda^\frac{1-d}{d} \Big[ A( \lambda I - P )^{-1}P - ( \lambda I - P )^{-1}PA \Big]u \, d\lambda
=0,
\end{align*}
and also 
\begin{align*}
[A, e^{-\epsilon P^\frac{1}{d}}]u
&=
A\Big(\frac{1}{2\pi i} \int_{\Gamma_R'} e^{-\epsilon \lambda^{\frac{1}{d}}} ( \lambda I - P )^{-1} u \, d\lambda \Big)
- e^{-\epsilon P^\frac{1}{d}}Au + 0 \\
&=
\frac{1}{2\pi i}\int_{\Gamma_R'} e^{-\epsilon \lambda^{\frac{1}{d}}} \Big[ A( \lambda I - P )^{-1} - ( \lambda I - P )^{-1} A \Big]u \, d\lambda
=0.
\end{align*}
\end{proof}

In the above lemma, $I+P$ is always invertible, $(I+P)^s$ is defined for any $s\in \mathbb{R}$. 
Of course, in that case, $A$ also commutes with any such power:
\begin{align*}
[A, (I+P)^s ] = 0,
\end{align*}
where the argument is similar to the above calculations.\\

Assume that $P$ and $\epsilon_0'\leq \epsilon_0$ are chosen so  Conjecture~\ref{conj:BoutetDeMonvelConjecture} holds, and $\epsilon\in (0,\epsilon'_0)$. 
The operator $P$ satisfies the conditions of Lemma~\ref{lmm:PowersExponentsCommutator}, so the calculus is applicable.
This is possible in at least one non-trivial case, $P=-\Delta_g$, by Theorem~\ref{thm:BoutetDeMonvelTheorem}.

\newpage
Let $h$ be a K\"ahler metric on the closure of $M^\Phi_\epsilon$, and let $s\in \mathbb{R}$ be arbitrary.

\begin{definition}
If $f\in H^{s-\frac{n+1}{4}}(M)$, $\mathcal{P}_\epsilon f$ is the extension of $e^{-\epsilon P^\frac{1}{d}}f$ to $M_\epsilon^\Phi$.
\end{definition}

Suppose $K_\epsilon$ is the Poisson operator solving the Dirichlet problem for $\Delta_h$ on $\overline{M_\epsilon^\Phi}$. As before, by Theorem~\ref{thm:PoissonSolutionOperator} it is a linear homeomorphism
\begin{align*}
K_\epsilon :H^{s-\frac{1}{2}}(\partial M_\epsilon^\Phi)  \to \bar{H}^s (M_\epsilon^\Phi) \cap \ker(\Delta_h |_{\bar{H}^s (M_\epsilon^\Phi)}),
\end{align*}
and $\mathcal{P}_\epsilon= K_\epsilon S_\epsilon$ realizes a well-defined linear homeomorphism 
\begin{align*}
\mathcal{P}_\epsilon : H^{s-\frac{n+1}{4}}(M) \to K_\epsilon \mathcal{O}^{s-\frac{1}{2}}(\partial M_\epsilon^\Phi).
\end{align*}

\begin{definition}
\begin{align*}
HH^s(M_\epsilon^\Phi) = K_\epsilon \mathcal{O}^{s-\frac{1}{2}}(\partial M_\epsilon^\Phi).
\end{align*}
\end{definition}

Note that $I+P\in \PsiPHG^d(M)$ is invertible, the inverse always belongs to $\PsiPHG^{-d}(M)$. By Lemma~\ref{lmm:PoissonImage} each space is a Hilbert space of holomorphic functions with norm
\begin{align*}
||u||_{HH^s(M^\Phi_\epsilon)}=||\mathcal{P}_\epsilon (I +P)^\frac{s}{d} \mathcal{P}_\epsilon^{-1} u||_{L^2(M_\epsilon^\Phi)}
\quad
\textnormal{if}
\quad
u \in HH^s(M^\Phi_\epsilon),
\end{align*}
and, automatically, $\mathcal{P}_\epsilon$ becomes a bounded isomorphism onto this space.

\begin{theorem} \label{thm:PAlgebraExtendible}
Let $d'\in \mathbb{R}$ and let $A\in \Psi^{d'}(M)$ commute with $P$, $[A,P]=0$. Then, if $s\in \mathbb{R}$, the diagram commutes, and consists of bounded operators:
\[
\begin{tikzcd}[row sep=large, column sep = huge]
HH^s(M_\epsilon^\Phi) \arrow[r, "\mathcal{P}_\epsilon A\mathcal{P}_\epsilon^{-1}"] \arrow[d, "\mathcal{R}_\epsilon"] & HH^{s-d'}(M_\epsilon^\Phi) \arrow[d, "\mathcal{R}_\epsilon"]  \\
H^{s}(M) \arrow[r, "A"] & H^{s-d'}(M)   
\end{tikzcd}
\]
\end{theorem}
\begin{proof}
By construction, we have
\begin{align*}
\mathcal{R}_\epsilon \mathcal{P}_\epsilon|_{H^{s-\frac{n+1}{4}}(M)} = e^{-\epsilon P^\frac{1}{d}}|_{H^{s-\frac{n+1}{4}}(M)}
\quad
\textnormal{for any}
\quad
s\in \mathbb{R},
\end{align*}
and so by Lemma~\ref{lmm:PowersExponentsCommutator}, if $u\in C^\infty(M)$ this implies
\begin{align*}
\mathcal{R}_\epsilon \mathcal{P}_\epsilon A u
=
e^{-\epsilon P^\frac{1}{d}} A u
=
A e^{-\epsilon P^\frac{1}{d}} u
=  
A \mathcal{R}_\epsilon \mathcal{P}_\epsilon u,
\end{align*}
which then extends by continuity to all $u\in H^s(M)$. Therefore the diagram commutes.
In particular, if $u\in HH^{s}(M_\epsilon)$, we have
\begin{align*}
(I+P)^{\frac{s}{d}} \mathcal{R}_\epsilon u = \mathcal{R}_\epsilon \mathcal{P}_\epsilon (I+P)^{\frac{s}{d}} \mathcal{P}_\epsilon^{-1} u,
\end{align*}
and so, by the analogue of \cite[pp.7, Inequalities (2.6)]{Stenzel2015} for $L^2(M^\Phi_\epsilon)$, we get
\begin{align*}
||\mathcal{R}_\epsilon u||_{H^{s}(M)} &\leq ||(I+P)^{-\frac{s}{d}}||_{B(H^{s}(M),L^2(M))} ||(I+P)^{\frac{s}{d}} \mathcal{R}_\epsilon u||_{L^2(M)} \\
&\leq
||(I+P)^{-\frac{s}{d}}||_{B(H^{s}(M),L^2(M))} ||\mathcal{R}_\epsilon||_{B(L^2(M^\Phi_\epsilon),L^2(M))} ||u||_{HH^{s}(M_\epsilon^\Phi)},
\end{align*}
which shows that the $\mathcal{R}_\epsilon$ are bounded in the diagram
\end{proof}

\newpage
The fact that the conjecture holds for $P=-\Delta_g$ leads to interesting corollaries. Those operators commuting with $\Delta_g$ always preserve $\epsilon$-extendible functions:

\begin{corollary} \label{cor:LaplacianAlgebraExtendible}
Let $d'\in \mathbb{R}$ and let $A\in \Psi^{d'}(M)$ commute with $\Delta_g$, $[A,\Delta_g]=0$. 
Then, if $s\in \mathbb{R}$, the diagram commutes, and consists of bounded operators:
\[
\begin{tikzcd}[row sep=large, column sep = huge]
HH^s(M_\epsilon) \arrow[r, "\mathcal{P}_\epsilon A\mathcal{P}_\epsilon^{-1}"] \arrow[d, "\mathcal{R}_\epsilon"] & HH^{s-d'}(M_\epsilon) \arrow[d, "\mathcal{R}_\epsilon"]  \\
H^{s}(M) \arrow[r, "A"] & H^{s-d'}(M)   
\end{tikzcd}
\]
\end{corollary}

\begin{proposition} \label{prop:LaplacianCommutatingOperators}
Let $\{\lambda_k\}_{k=1}^\infty$ be the eigenvalues of $-\Delta_g$, with multiplicities, and $\{\psi_k\}_{k=1}^\infty$ an orthonormal eigenbasis enumerated according to these eigenvalues.
Suppose $ \sup_{k\in \mathbb{N}} \langle\lambda_k\rangle^{-d'}|\eta_k| < \infty$ for some $d'\in \mathbb{R}$, and put
\begin{align*}
A : C^\infty(M) \to C^\infty(M) : u \mapsto \sum_{k=1}^\infty \eta_k (u, \psi_k)_{L^2(M)} \psi_k.
\end{align*}
Then $A$ is continuous, 
$[A, \Delta_g]=0$, and the conclusion of Corollary~\ref{cor:LaplacianAlgebraExtendible} holds.
\end{proposition}
\begin{proof}
Taking any $s\in \mathbb{R}$ and $u\in H^s(M)$, we have
\begin{align*}
||Au||_{H^{s-d'}(M)}^2 &=
\sum_{k=1}^\infty \langle \lambda_k \rangle^{2(s-d')}|\eta_k|^2 | ( u, \psi_k )_{L^2(M)} |^2 \\
&\leq 
C \sum_{k=1}^\infty \Big| \Big( (I-\Delta_g)^{-\frac{s}{2}} u, \psi_k \Big)_{L^2(M)} \Big|^2 
=C ||u||_{H^s(M)}^2,
\end{align*}
and so by the Sobolev embedding theorem, $A$ is continuous from $C^\infty(M)$ into itself. 
It is not necessary that $A \in \Psi(M)$, we can use  
Lemma~\ref{lmm:PowersExponentsCommutator} just as in Theorem~\ref{thm:PAlgebraExtendible}.
Repeating the arguments lead to the conclusions of Corollary~\ref{cor:LaplacianAlgebraExtendible}
\end{proof}

\begin{corollary} \label{cor:laplacianparametrixresult}
$-\Delta_g$ admits a parametrix for which the diagram commutes. 
\end{corollary}
\begin{proof}
Put $\eta_k = \frac{1}{\lambda_k}$ when $\lambda_k\neq 0$ and zero otherwise.
\end{proof}

Thus, the commutant of $\Delta_g$ consists of operators with the $\epsilon$-extension property. 
In particular, this includes any operator in the functional calculus of $\Delta_g$.\\

The most interesting implication of the Conjecture~\ref{conj:BoutetDeMonvelConjecture} is a generalization of this. Let $d' > 0$ and 
assume that $A\in \PsiPHG^{d'}(M)$ is analytic, elliptic and formally normal. In that case, let $a$ be the classical principal symbol of $A$, and put
\begin{align*}
P = A^*A\in \PsiPHG^{d}(M)
\quad
\textnormal{with}
\quad
d=2d',
\end{align*}
which then has the unique classical principal symbol
\begin{align*}
p=|a|^2
\quad
\textnormal{and}
\quad
p|_{T^*M \setminus 0} >0.
\end{align*}
This $P$ is elliptic, formally self-adjoint, and can not have any negative eigenvalues. But, we still need a convexity assumption to satisfy the conditions of Conjecture~\ref{conj:BoutetDeMonvelConjecture}, and so we
require of $a$ that $\{ \xi \in T^*_x M \, | \, |a(x,\xi)|\leq 1 \}$ is strictly convex for any $x\in M$.

\newpage
\begin{theorem} \label{thm:parametrixresult}
$A$ admits a parametrix $B\in \PsiPHG^{-d'}(M)$ such that $[A^*A, B]=0$. Consequently, if Conjecture~\ref{conj:BoutetDeMonvelConjecture} holds for $P=A^*A$, the diagram commutes:
\[
\begin{tikzcd}[row sep=large, column sep = huge]
HH^s(M_\epsilon^\Phi) \arrow[r, "\mathcal{P}_\epsilon B\mathcal{P}_\epsilon^{-1}"] \arrow[d, "\mathcal{R}_\epsilon"] & HH^{s+d'}(M_\epsilon^\Phi) \arrow[d, "\mathcal{R}_\epsilon"]  \\
H^{s}(M) \arrow[r, "B"] & H^{s+d'}(M)   
\end{tikzcd}
\]
\end{theorem}
\begin{proof}
If $A$ is invertible, $A^{-1}$ commutes with $P$ on $C^\infty(M)$, and so we are done. Assume that $A$ is not invertible.
Then $0$ is an eigenvalue of multiplicity $\dim \ker P < \infty$, and
$\sigma(P)\subset [0,\infty)$ consists of eigenvalues, $0=\lambda_1 \leq \lambda_2 \cdots$, counted with multiplicity.
By \cite[Theorem 9.3]{Shubin2001}, $\int_\Gamma \frac{1}{\lambda} (\lambda I - P)^{-1}\, d\lambda$ converges in $B(H^N(M))$ for all $N\in \mathbb{N}$, where $\Gamma$ is a circle inside the resolvent set of $P$, enclosing $0$ but no other eigenvalues.
Then by the Sobolev embedding theorem, we
get a continuous operator
\begin{align*}
B : C^\infty(M) \to C^\infty(M) : u \mapsto  -\frac{1}{2\pi i} \int_\Gamma \frac{1}{\lambda} (\lambda I - P)^{-1}A^*u \, d\lambda.
\end{align*}
Take $u\in C^\infty(M)$, and use the formal normality of $A$ to write
\begin{align*}
[A^*A, B]u &= -\frac{1}{2\pi i} P \Big( \int_\Gamma \frac{1}{\lambda} (\lambda I - P)^{-1}A^*u \, d\lambda \Big) - BP u \\
&=
-\frac{1}{2\pi i}  \int_\Gamma \frac{1}{\lambda} \Big[ P(\lambda I - P)^{-1}A^* - (\lambda I - P)^{-1}A^*P \Big] u \, d\lambda = 0,
\end{align*}
and, if $\{ \psi_k \}_{k=1}^\infty$ is an orthonormal eigenbasis corresponding to $\{\lambda_k\}_{k=1}^\infty$, we have 
\begin{align*}
AB u = BA u
&= -\frac{1}{2\pi i} \int_\Gamma \frac{1}{\lambda} (\lambda I - P)^{-1}P u \, d\lambda \\
&=
\Big[\frac{1}{2\pi i} \int_\Gamma \frac{1}{\lambda}  \, d\lambda \Big] u - \frac{1}{2\pi i} \int_\Gamma (\lambda I - P)^{-1} \Big[ \sum^\infty_{k=1} (u, \psi_k)_{L^2(M)} \psi_k \Big] \, d\lambda  \\
&=
Iu - \sum^\infty_{k=1} (u,\psi_k)_{L^2(M)} \Big[ \frac{1}{2\pi i} \int_\Gamma \frac{1}{\lambda - \lambda_k}  \, d\lambda \Big] \psi_k \\
&=
Iu - \sum^{\dim \ker P}_{k=1} (u,\psi_k)_{L^2(M)} \psi_k,
\end{align*}
where the remainder is simply the projection onto the finite-dimensional kernel of $A$. 
So $B$ can only differ from a parametrix $B'\in \PsiPHG^{-d'}(M)$ of $A$ by a smoothing operator:
\begin{align*}
(B-B')u 
&= B(I-AB')u-(I-BA)B'u \\
 &= \int_M \Big( B [R(\cdot,y)](x)   - \sum^{\dim \ker P}_{k=1} \psi_k(x) ((B')^*\psi_k)(y) \Big) u(y) \omega_0(y),
\end{align*}
where $\omega_0$ is the 1-density relative to which $L^2(M)$ is defined, $R$ the kernel of $I-AB'$. 
It follows that $B\in \PsiPHG^{-d'}(M)$, and that it is a parametrix for $A$. 
\end{proof}

\newpage
An analogue of \cite[pp. 6, Theorem 1]{Stenzel2015} holds for the holomorphic Sobolev spaces. The proof is very similar, but we generalize it, and fill in some details.


\begin{theorem} Suppose that
Conjecture~\ref{conj:BoutetDeMonvelConjecture} holds for an operator $P$ as in \ref{conj:BoutetDeMonvelConjecture}, and that the classical principal symbol of $S^*_\epsilon S_\epsilon \in \PsiPHG^{-\frac{n-1}{2}}(M)$ is positive on $T^* M \setminus 0$. 
Then, given any $s\in \mathbb{R}$, the following holds:
\medskip 
\begin{enumerate}
\item There is a well-defined, positive operator
\begin{align*}
\mathcal{Q}_\epsilon = ( (I+P)^\frac{s}{d} S_\epsilon^* (K_\epsilon^* K_\epsilon) S_\epsilon (I+P)^\frac{s}{d} )^{-\frac{1}{2}} (I - \Delta_g)^\frac{s}{2}
\in \PsiPHG^{\frac{n+1}{4}}(M),
\end{align*}

\medskip 

\item It unitarizes $\mathcal{P}_\epsilon$, in that sense that 
\begin{align*}
||\mathcal{P}_\epsilon \mathcal{Q}_\epsilon u ||_{HH^s(M^\Phi_\epsilon)} = ||u||_{H^s(M)}
\quad
\textnormal{for any}
\quad
u\in H^s(M).
\end{align*}

\end{enumerate}
\medskip 
Here $K_\epsilon^*$ is understood to be the formal $L^2$-adjoint of $K_\epsilon $ obtained from Theorem~\ref{thm:PoissonSolutionOperator}, and $H^s(M)$ is normed so that $(I - \Delta_g)^\frac{s}{2} : H^s(M) \to L^2(M)$ is an isometry.
\end{theorem}

\begin{proof}
Using Theorems~\ref{thm:PoissonSolutionOperator} and \ref{thm:BoutetDeMonvelTheorem} and the complex Egorov theorem, we get 
\begin{align*}
S_\epsilon^* (K^*_\epsilon K_\epsilon) S_\epsilon \in \PsiPHG^{-\frac{n+1}{2}}(M),
\end{align*}
which is then elliptic, and its classical principal symbol is strictly positive on $T^* M \setminus 0$.
It has no negative eigenvalues. We show that it is also invertible.\\

Regarding $\mathcal{P}_\epsilon$ as a bounded map into $L^2(M_\epsilon^\Phi)$, we form its Hilbert adjoint $\mathcal{P}_\epsilon^*$.
Then $\mathcal{P}_\epsilon^* \mathcal{P}_\epsilon : H^{-\frac{n+1}{4}}(M) \to H^{-\frac{n+1}{4}}(M) $ is injective with dense range, and
\begin{align*}
\mathcal{P}_\epsilon^* \mathcal{P}_\epsilon|_{C^\infty(M)}
=
(I-\Delta_g )^{\frac{n+1}{4}} S_\epsilon^* (K_\epsilon^* K_\epsilon) S_\epsilon \in \PsiPHG^0(M).
\end{align*}
But this shows that $\mathcal{P}_\epsilon^* \mathcal{P}_\epsilon$ is a Fredholm operator on $H^{-\frac{n+1}{4}}(M)$, so has closed range.
Therefore, we have a bijective and continuous, hence invertible, operator
\begin{align*}
S_\epsilon^* (K^*_\epsilon K_\epsilon) S_\epsilon : C^\infty(M) \to C^\infty(M),
\end{align*}
and the inverse is again a pseudo-differential operator by \cite[pp. 70, Theorem 8.2]{Shubin2001}.
Applying \cite[pp. 90-91, Theorem 10.1]{Shubin2001} and \cite[pp. 93, Proposition 10.3]{Shubin2001}, $\mathcal{Q}_\epsilon$ exists.
It follows that $\mathcal{Q}_\epsilon$ is well-defined, elliptic, invertible, and also formally self-adjoint. 
Taking any $u\in C^\infty(M)$, we see that
\begin{align*}
\big( \mathcal{Q}_\epsilon u , u \big)_{L^2(M_\epsilon^\Phi)}
=
\big( \mathcal{Q}_\epsilon^\frac{1}{2} u , \mathcal{Q}_\epsilon^\frac{1}{2} u \big)_{L^2(M_\epsilon^\Phi)} > 0
\quad
\textnormal{if}
\quad
u \neq 0,
\end{align*}
and
\begin{align*}
||(K_\epsilon S_\epsilon)  u ||^2_{HH^s(M_\epsilon^\Phi)}
&=
\big( (K_\epsilon S_\epsilon)(I+P)^\frac{s}{d} u , (K_\epsilon S_\epsilon)(I+P)^\frac{s}{d} u \big)_{L^2(M_\epsilon^\Phi)} \\
&=
\big( (I+P)^\frac{s}{d}(K_\epsilon S_\epsilon)^*(K_\epsilon S_\epsilon)(I+P)^\frac{s}{d} u ,  u \big)_{L^2(M)}  \\
&=
\big( (I-\Delta_g)^{\frac{s}{2}} \mathcal{Q}_\epsilon^{-1} u , (I-\Delta_g)^{\frac{s}{2}} \mathcal{Q}_\epsilon^{-1} u \big)_{L^2(M)} \\
&= || \mathcal{Q}_\epsilon^{-1} u||^2_{H^s(M)}.
\end{align*}
\end{proof}

\newpage
\subsection{On a Compact Lie Group $G$}
Let $\{G_\epsilon\}_{\epsilon>0}$ be the Grauert tubes of $G$.
It is possible to form a subalgebra of $\Psi(G)$ with the desired properties for any $\epsilon>0$. 
That is, operators acting as bounded maps between the holomorphic Sobolev spaces, and 
with an analogous notion of, and conditions for, ellipticity within the subalgebra.
We will show that the subalgebra is non-trivial. 
Let $d\in \mathbb{R}$.

\begin{definition} \label{def:Sepsilon}
Define $S^d_\epsilon$ to be those $p\in S^d$ satisfying:
\medskip
\begin{enumerate}

\item The mapping $G \to \mathfrak{gl}(d_\xi, \mathbb{C}) : x \mapsto p(x,\xi)$ is real-analytic for every $[\xi]\in \widehat{G}$. Each of these extend holomorphically to $G_\epsilon$.

\medskip

\item Write $p_Y(x,\xi)=p(\exp(iY)x,\xi)$ for every $(x,[\xi])\in G\times \widehat{G}$ when $|Y|_\mathfrak{g}<\epsilon$. Then $\{ p_Y \}_{|Y|_\mathfrak{g}<\epsilon}$ is a bounded subset of $S^d$. 

\end{enumerate}
\medskip
\end{definition}

Take $p\in S^d_\epsilon$. We employ the notation $p_Y$ as it appears above in the second point.
If $u\in C^\omega(G)$ extends to $G_\epsilon$, then since $G$ is totally real in $G_\mathbb{C}$, we have
\begin{align*}
\Op(p)u(\exp(iY)x) 
&= 
\sum_{[\xi] \in \widehat{G}} d_\xi \Tr(\xi(\exp(iY)x)p_Y(x,\xi)\mathcal{F}u(\xi)) \\
&= 
\sum_{[\xi] \in \widehat{G}} d_\xi \Tr\Big(\xi(x)p_Y(x,\xi)\int_G u(\exp(iY)y) \xi(y)^* \, dy \Big), 
\end{align*}
where the sum is then uniformly absolutely convergent on $x\in G $ and $|Y|\leq r < \epsilon$.
This gives us the following simple result.

\begin{proposition}
If $u\in C^\omega(G)$ extends holomorphically to $G_\epsilon$, so does $\Op(p)u$.
\end{proposition}

\begin{theorem}
Let $d_1, d_2\in \mathbb{R}$. If $p\in S^{d_1}_\epsilon$ and $q\in S^{d_2}_\epsilon$, then $p\odot q \in S^{d_1+d_2}_\epsilon$.
\end{theorem}
\begin{proof}
If $[\eta]\in \widehat{G}$ is fixed, then $p(x,\eta)$ extends holomorphically in $x\in G$ to $z\in G_\epsilon$. This is clear from the formula
\begin{align*}
(p\odot q)(z,\eta)
&=
\sum_{[\xi] \in \widehat{G}}
d_\xi \int_G  \Tr\Big(\xi(y) p(z ,\xi) \Big) \eta(y)^* q(zy^{-1} ,\eta) \, dy,
\end{align*}
and, if $|Y|_\mathfrak{g}<\epsilon$, it also shows that
\begin{align*}
(p\odot q)_Y (x,\eta)
&=
(p\odot q)(\exp(iY)x,\eta) \\
&=
(p_Y \odot q_Y) (x,\eta).
\end{align*}
Then, because $\{p_Y\}_{|Y|_\mathfrak{g}<\epsilon}\subset S^{d_1}$ and  $\{p_Y\}_{|Y|_\mathfrak{g}<\epsilon}\subset S^{d_2}$ are bounded sets by definition, and the symbolic product is continuous, $\{ (p\odot q)_Y \}_{|Y|_\mathfrak{g}<\epsilon} $ is a bounded set in $S^{d_1+d_2}$.
This shows that $p \odot q \in S^{d_1+d_2}_\epsilon$.
\end{proof}

It therefore follows that $\Psi_\epsilon(G)=\cup_{d\in \mathbb{R}}\Op\, S^d_\epsilon$ is actually a subalgebra of $\Psi(G)$. 
The holomorphic extension of $\Op(p)u$ can also be expressed in terms of $\{\Op(p_Y)\}_{|Y|_\mathfrak{g}<\epsilon}$. In fact, with $x$ and $Y$ as above, and $L_{(\cdot)}$ the left regular representation, we have
\begin{align*}
\Op(p)u(z) = \Op(p_Y) (L_{\exp(-iY)}u)(x)
\quad
\textnormal{if}
\quad
z=\exp(iY)x \in G_\epsilon,
\end{align*}
which by itself is \textit{not} a pseudo-differential operator.

\newpage
Next, we get the holomorphic mapping property using the Cartan decomposition. It mirrors precisely the standard Sobolev mapping property of $\Psi(G)$.

\begin{theorem} \label{thm:LieGroupSobolevMirror}
If $p\in S^d_\epsilon$ and $s\in \mathbb{R}$, then
$\Op(p) $ admits a unique $\epsilon$-extension.
That is, the diagram commutes, and consists of bounded operators:
\[
\begin{tikzcd}[row sep=large, column sep = huge]
HH^s(G_\epsilon) \arrow[r, "\widetilde{\Op(p)}"] \arrow[d, "\mathcal{R}_\epsilon"] & HH^{s-d}(G_\epsilon) \arrow[d, "\mathcal{R}_\epsilon"]  \\
H^{s}(G) \arrow[r, "\Op(p)"] & H^{s-d}(G)   
\end{tikzcd}
\]
\end{theorem}
\begin{proof}
Because $(I-\Delta_G)^\frac{s}{2}\in \Op S^s_\epsilon$ for all $s\in \mathbb{R}$ it is enough to prove $s=d=0$. Using the boundedness of $\{p_Y\}_{|Y|_\mathfrak{g}<\epsilon}$, Lemma~\ref{lmm:LieGroupCalderonVaillancourt}, and \cite[Lemma 5]{Hall1997}, we get
\begin{align*}
\int_{G_\epsilon} |\Op(p)u(z)|^2 \, dz
&=
\int_{|Y|_\mathfrak{g}<\epsilon} \Big[ \int_G |[\Op(p)u](\exp(iY)x)|^2 \, dx  \Big] \, \frac{dY}{\Theta(Y)^2} \\
&=
\int_{|Y|_\mathfrak{g}<\epsilon} \Big[ \int_G |[\Op(p_Y)L_{\exp(-iY)}u](x)|^2 \, dx  \Big] \, \frac{dY}{\Theta(Y)^2} \\
&\leq
C
\int_{|Y|_\mathfrak{g}<\epsilon} \Big[ \int_G |u(\exp(iY)x)|^2 \, dx  \Big] \, \frac{dY}{\Theta(Y)^2} \\
&= 
C \int_{G_\epsilon} |u(z)|^2 \, dz,
\end{align*}
where $C>0$ is the constant in the uniform bound
\begin{align*}
|| \Op(p_Y)L_{\exp(-iY)}u ||_{L^2(G)}^2 \leq C || L_{\exp(-iY)}u ||_{L^2(G)}^2,
\end{align*}
and $Y \mapsto \Theta(Y)^{-2}$ is the $\Ad(G)$-invariant Jacobian from \cite[Lemma 5]{Hall1997}.
\end{proof}

\begin{remark}
Take any two symbols $p\in S^{d_1}_\epsilon$ and $q\in S^{d_2}_\epsilon$ for some $d_1,d_2\in \mathbb{R}$. Then $(p\odot q)_Y = p_Y \odot q_Y$ so $\{ p_Y\odot q_Y - p_Y q_Y \}_{|Y|_\mathfrak{g}<\epsilon} \subset S^{d_1+d_2-1}$ is a bounded set, and we must have
\begin{align*}
p\odot q - pq \in S^{d_1+d_2-1}_\epsilon.
\end{align*}
It does not appear possible to include more terms, as these involve the $\delta_x$ operator, which contains a cutoff, so the higher order terms may not extend holomorphically. Nevertheless, we only need this to hold for the principal term.
\end{remark}

\begin{proposition}
The space $S^d_\epsilon$ contains non-trivial symbols for each $d\in \mathbb{R}$.
\end{proposition}
\begin{proof}
The weighted Bargmann space $HL^2(G_\mathbb{C}, \nu_t) $ of \cite{Hall1994} is infinite-dimensional.
Hence so is the space of the restrictions of these functions to $G_\epsilon$, and so also $\mathcal{O}(G_{\epsilon})$.
Let $S^d_{\textnormal{inv}}$ be the $d\in \mathbb{R}$ symbols of bi-invariant operators on $G$ (depending only on $[\xi]$). Then, we have that
\begin{align*}
\mathcal{O}(G_{\epsilon'})|_{G_\epsilon} \otimes S^d_{\textnormal{inv}}\subset S^d_\epsilon
\quad
\textnormal{when}
\quad
\epsilon'>\epsilon,
\end{align*}
and this space is infinite-dimensional, the symbols depend on both $z$ and $[\xi]$.
\end{proof}

\newpage
The spaces $\Op \, S^k_\epsilon$ contain the real-analytic differential operators of degree $k\in \mathbb{N}$. This is easy to prove by writing out their matrix-symbols. Put
\begin{align*}
\mathfrak{g}_\epsilon = \{ Y\in \mathfrak{g} \, | \, |Y| <\epsilon \}.
\end{align*}
In the following, $(X_1, \cdots, X_n)$ is an ordered basis of left-invariant vector fields for $\mathfrak{g}$, and we write $X^\beta = X^{\beta_1}_1 \circ \cdots \circ X^{\beta_n}_n$ for any $\beta = (\beta_1, \cdots, \beta_n)\in \mathbb{N}_0^n$. 

\begin{proposition} \label{prop:realandifflie}
If $P\in \Diff^k(G)$ is real-analytic, $P\in \Op \, S^k_\epsilon$ for some $\epsilon>0$.
\end{proposition}
\begin{proof}
In any sufficiently small chart $U$ of $G$, we have
\begin{align*}
P|_U = \sum_{|\beta|\leq k} g_\beta X^\beta,
\end{align*}
where $g_\beta \in C^\omega(U)$ extend holomorphically into $\exp(i \mathfrak{g}_\epsilon) U$ for $\epsilon>0$ depending on $U$.
It follows that if $z= \exp(iY)x \in \exp(i \mathfrak{g}_\epsilon) U$, we have
\begin{align*}
\xi(z)^{-1}P\xi(z) 
&=
\xi(x)^{-1}\xi(\exp(iY))^{-1}\sum_{|\beta|\leq k} g_\alpha(z) (X^\beta \xi)(\exp(iY)x) \\
&=
\sum_{|\beta|\leq k} g_\beta (\exp(iY)x) \Big[ \xi(x)^{-1} X^\beta \xi(x) \Big]
\\ 
&= 
\sum_{|\beta|\leq k} g_\alpha(\exp(iY)x) p_\beta(\xi),
\end{align*}
and $(x,Y) \mapsto g_\alpha(\exp(iY)x)$ has all its $x$-derivatives uniformly bounded in $Y \in \mathfrak{g}_\epsilon$.
The symbols $p_\beta \in S^{|\beta|}$ are independent of all $(x,Y)\in G \times\mathfrak{g}_\epsilon$, and given by
\begin{align*}
p_\beta(\xi)=\xi(x)^{-1} X^\beta \xi(x)
=(X^\beta\xi)(e).
\end{align*}
Taking a finite cover of $G$ by such $U$,
we get $\epsilon>0$ so $\{p_Y\}_{Y \in \mathfrak{g}_\epsilon}$ is bounded in $S^k$.
Thus we have $P\in \Op \, S^d_\epsilon$ for some $\epsilon>0$.
\end{proof}

\begin{corollary}
If $P$ above is left-invariant, then $P\in \Op\, S^k_\epsilon$ for all $\epsilon>0$.
\end{corollary}
\begin{proof}
In this case, the coefficients $g_\alpha$ are globally defined and constant.
\end{proof}

To characterize ellipticity for these symbols, we need to sum them asymptotically.
It will be necessary to use the following, similar to \cite[pp. 20, Lemma 7.2]{RuzhanskyWirth2014}.

\begin{lemma} \label{lmm:newsymbtransf}
Let $d\in \mathbb{R}$ and $p\in S^d_\epsilon$, and put
\begin{align*}
R(z)=\mathcal{F}^{-1}_{\xi}[p(z,\xi)]\in \mathcal{D}'(G)
\quad
\textnormal{for any}
\quad
z\in G_\epsilon,
\end{align*}
where $(x,Y) \mapsto p(\exp(iY)x,\xi)$ has been extended by continuity to $x\in G$ and $|Y|\leq \epsilon$.
It is alternatively viewed as a map
\begin{align*}
G\times \overline{\mathfrak{g}_\epsilon} \to \mathcal{D}'(G) : (x,Y) \mapsto R(\exp(iY)x)=R(x,Y).
\end{align*}
Then, we have
\begin{align*}
R \in C^\infty(G\times \mathfrak{g}_\epsilon, H^{-d-\lceil \frac{n}{2} \rceil}(G)),
\end{align*}
and the $x$-derivatives of $R$ are uniformly bounded in $Y$ up to $\overline{\mathfrak{g}_\epsilon}$.
\end{lemma}

\newpage
\begin{proof}
Moving derivatives under $\mathcal{F}^{-1}_\xi$, we can differentiate $R$ in the weak sense. The $X^\beta_x$ derivatives are just
\begin{align*}
X^\beta_x R(x,Y) = \mathcal{F}^{-1}_{\xi}[X_x^\beta p_Y ( x,\xi)]
\quad
\textnormal{in}
\quad
C^\infty(G \times \mathfrak{g}_\epsilon, \mathcal{D}'(G)).
\end{align*}
Using the Parseval identity, if $x\in G$, we have
\begin{align*}
||R(x,Y)||_{H^{-d-\lceil \frac{n}{2} \rceil}(G)} &=
\sum_{[\xi]\in \widehat{G}} d_\xi \langle \xi \rangle^{-2d-2\lceil \frac{n}{2} \rceil} \Tr \Big(  p_Y(x,\xi)^* p_Y (x,\xi) \Big) \\
&\leq
\sum_{[\xi]\in \widehat{G}} d_\xi^2 \langle \xi \rangle^{-2\lceil \frac{n}{2} \rceil} \Big[ \langle \xi \rangle^{-d} || p_Y(x,\xi) || \Big]^2 < \infty,
\end{align*}
and, similarly, for any weak $x$-derivative of $x \mapsto R(x,Y)$, we find it in $H^{-d-\lceil \frac{n}{2} \rceil}(G)$. So $ (x,Y) \mapsto R(x,Y)\in H^{-d-\lceil \frac{n}{2} \rceil}(G)$ has weak $x$-derivatives uniformly bounded in $Y$.
It remains to show smoothness in norm, the derivatives will agree with the weak ones.
In fact, even more is true. The function $z \mapsto R(z)$ is strongly holomorphic. \\

It is sufficient to prove weak holomorphy due to
\cite[pp. 82-83, Theorem 3.31]{RudinFunctionalAnalysis}.
At any $z_0\in G_\epsilon$ write $z=z_0\exp(W)$ for $W$ in a small polydisc about $0$ in $\mathfrak{g} \otimes_\mathbb{R} \mathbb{C}$.
Take $v\in H^{d+\lceil \frac{n}{2} \rceil}(G)$ and a closed piecewise $C^1$ curve $\gamma$ in the $j$-th disc, and write
\begin{align*}
\int_\gamma \langle v, R(z) \rangle \, dW_j
&=
\int_\gamma \langle v, \mathcal{F}^{-1}_{\xi}[ p ( z,\xi)]  \rangle \, dW_j \\
&=
\int_\gamma \Big[ \int_G (I-\Delta_G)^{\frac{d+\lceil \frac{n}{2} \rceil}{2}} v(x) \mathcal{F}^{-1}_{\xi}[\langle \xi\rangle^{-\frac{d+\lceil \frac{n}{2} \rceil}{2}} p ( z,\xi)]  \, dx \Big] \, dW_j \\
&=
\int_\gamma \Big[ \sum_{[\xi] \in \widehat{G}} d_\xi \Tr\Big( \mathcal{F}\overline{v}(\xi)^* p(z,\xi) \Big) \Big] \, dW_j \\
&=
\sum_{[\xi] \in \widehat{G}} d_\xi \Tr\Big( \mathcal{F}\overline{v}(\xi)^* \int_\gamma p(z,\xi) \, dW_j \Big) 
= 0.
\end{align*}
It follows from Morera's theorem and Osgood's lemma that
\begin{align*}
R \in \mathcal{O}(G_\epsilon, H^{-d-\lceil \frac{n}{2} \rceil}(G)),
\end{align*}
and $G \times \mathfrak{g}_\epsilon \to H^{-d-\lceil \frac{n}{2} \rceil}(G) : (x,Y) \mapsto R(x,Y)$ is then norm-differentiable a fortiori. The $x$-derivatives are uniformly bounded in $Y$ up to $\overline{\mathfrak{g}_\epsilon}$, since the weak ones are. 
\end{proof}

\begin{lemma} \label{lmm:asymptoticsumcomplex}
Let $\{p_j\}_{j=0}^\infty$ be a sequence with $p_j\in S^{d_j}_\epsilon$
and $d_j \searrow -\infty$ as $j\to \infty$. Then we can construct $p\in S^{d_0}_\epsilon$ with $p \sim \sum_{j=0}^\infty p_j$ in the sense that
\begin{align*}
p- \sum^{k-1}_{j=0} p_j \in S^{d_k}_\epsilon
\quad
\textnormal{for each}
\quad
k\in \mathbb{N}.
\end{align*}
\end{lemma}

\begin{proof}
In this proof, we let $\beta\in \mathbb{N}^n_0$ and $\alpha \in \mathbb{N}^{m\times m}_0$ be arbitrary multi-indices. Take an approximate identity $\{ \psi_t \}_{t>0} \subset C^\infty(G)$ for the convolution product $\ast$ on $G$. In this way, if $f\in H^s(G)$ for any $s\in \mathbb{R}$, we have
\begin{align*}
||\psi_{t}\ast f - f||_{H^s(G)} \to 0
\quad
\textnormal{as}
\quad
t \to 0.
\end{align*}

\newpage
Construct functions $R_j : G \times \overline{\mathfrak{g}_\epsilon} \to \mathcal{D}'(G)$ from $p_j \in S^{d_j}_\epsilon$ just as in Lemma~\ref{lmm:newsymbtransf}. Then $R_j \in C^\infty(G\times \mathfrak{g}_\epsilon, H^{-d_j-\lceil \frac{n}{2} \rceil}(G))$, and we can
pick $\{ t_j \}_{j=0}^\infty$ such that
\begin{align*}
\sup_{(x,Y)\in G\times   \overline{\mathfrak{g}_\epsilon}   }||X^\beta_x R_j(x,Y) - \psi_{t_j} \ast (X^\beta_x R_j(x,Y)) ||_{H^{-d_j-\lceil \frac{n}{2} \rceil}(G)} < \frac{1}{2^{j+1}}
\quad
\textnormal{when}
\quad
|\beta|\leq j.
\end{align*}
Take $N \in \mathbb{N}_0$ with $-d_j -\lceil \frac{n}{2} \rceil \geq 0$ for $j\geq N$, and observe that   
\begin{align*}
||p_j(z,\xi) (1-\mathcal{F}\psi_{t_j}(\xi))||^2
&\leq
||p_j(z,\cdot) (1-\mathcal{F}\psi_{t_j}(\cdot))||^2_{\mathcal{S}^0 (\widehat{G})} \\
&=
 ||R_j(z) - \psi_{t_j} \ast ( R_j(z)) ||_{L^2(G)}^2 \\
&\leq
||R_j(z) - \psi_{t_j} \ast ( R_j(z)) ||_{H^{-d_j-\lceil \frac{n}{2} \rceil}(G)}^2
< \frac{1}{4^{j+1}},
\end{align*}
where the estimates are uniform in $z=\exp(iY)x$ for $(x,Y)\in G \times \overline{\mathfrak{g}_\epsilon}$ and $[\xi]\in \widehat{G}$. Consequently, we have
\begin{align*}
\sum_{j= N}^\infty ||p_j(z,\xi) (1-\mathcal{F}\psi_{t_j}(\xi))|| < 1.
\end{align*}
It follows that $p\in S^{d_0}_\epsilon$ if we put 
\begin{align*}
p(z,\xi) = \sum_{j= 0}^\infty p_j(z,\xi) (1-\mathcal{F}\psi_{t_j}(\xi))
\quad
\textnormal{for any}
\quad
(z, [\xi]) \in G_\epsilon \times \widehat{G},
\end{align*}
which converges absolutely uniformly (in matrix norm) on $z\in G_\epsilon$ for each $[\xi]\in \widehat{G}$, and $x\mapsto p_Y(x,\xi)$ is differentiable, $X_x^\beta$ falls onto $x \mapsto p_j(\exp(iY)x,\xi)$ under the sum.
Let $r_N$ denote the sum starting at $j=N$. Take $k\in \mathbb{N}$, and write
\begin{align*}
p(z,\xi) - \sum^{k-1}_{j=0} p_j(z,\xi) = 
\sum^{N-1}_{j=k} p_j(z,\xi) 
- \sum^{N-1}_{j=k} p_j(z,\xi) \mathcal{F}\psi_{t_j}(\xi)
+ r_N(z,\xi),
\end{align*}
where the first term is in $S^{d_k}_\epsilon$, and the second belongs to $S^{-\infty}_\epsilon$ by Lemma~\ref{lmm:DifferenceOperatorLeibnizRule} for $\delta^\alpha_\xi$. 
It suffices that $r_N \in S^{d_k}_\epsilon$ for large $N\in \mathbb{N}_0$. But by \cite[Lemma 7.1]{RuzhanskyWirth2014}, we get
\begin{align*}
\langle \xi \rangle^{|\alpha|-d_k} || \delta^\alpha_\xi X^\beta_x (r_N )_Y(x,\xi) || 
&\leq 
C_{\alpha,\beta} \Big[\sup_{[\xi]\in \widehat{G}} \langle \xi \rangle^{|\alpha|-d_k} || X_x^\beta (r_N)_Y(x,\xi)|| \Big] \\
&\leq 
C_{\alpha,\beta} \Big[ \sum_{j = N}^\infty || X_x^\beta (p_j)_Y(x,\cdot) (1-\mathcal{F}\psi_{t_j}(\cdot)) ||_{H^{|\alpha|-d_k}(\widehat{G})} \Big] \\
&=
C_{\alpha,\beta} \Big[ \underbrace{\sum_{j = N}^\infty || X_x^\beta R_j(z) - \psi_{t_j} \ast X_x^\beta R_j(z) ||_{H^{|\alpha|-d_k}(G)} }_{\leq 1} \Big],
\end{align*}
where $N$ is chosen so that 
\begin{align*}
d_N + \lceil \frac{n}{2} \rceil < d_k -|\alpha|,
\end{align*}
and $C_{\alpha,\beta}>0$ is a constant depending only on $\alpha$ and $\beta$. 
\end{proof}

\newpage
The above means that we can sum symbols in the $\{S^d_\epsilon\}_{d\in \mathbb{R}}$ classes asymptotically. This allows us to get inverses in $S^{-d}_\epsilon$ mod $S^{-\infty}_\epsilon$ for elliptic symbols in $S^d_\epsilon$.

\begin{theorem} 
Let $q_0\in S^{-d}_\epsilon$.
The following holds:
\medskip
\begin{enumerate}

\item If $q_0p - 1 \in S^{-1}_\epsilon$, then there is a $q_L \in S^{-d}_\epsilon$ such that $q_L \odot p - 1 \in S^{-\infty}_\epsilon$. 

\medskip

\item If $p q_0  - 1 \in S^{-1}_\epsilon$, then there is a $q_R \in S^{-d}_\epsilon$ such that $p \odot q_R - 1 \in S^{-\infty}_\epsilon$. 

\end{enumerate}
\medskip
Finally, if both $q_L$ and $q_R$ as above exist, then $q_L - q_R \in S^{-\infty}_\epsilon$.
\end{theorem}

\begin{proof}
If both exist, then $q_L - q_R = q_L \odot ( 1- p \odot q_R ) - (1 - q_L \odot  p ) \odot q_R \in S^{-\infty}_\epsilon$.
Let $N\in \mathbb{N}$. In either case, Lemma~\ref{lmm:asymptoticsumcomplex} allows us to build the inverse $q$ as follows.\\


{\bf Left:} Put $r = 1- q_0 \odot p $. Define the symbol sequence $q_j = r^{\odot j}\odot q_0$ for $j\in \mathbb{N}_0$. Then put $q \sim \sum_{j=0}^\infty q_j $ with $q \in S^{-d}_\epsilon$, and write
\begin{align*}
S^{-N}_\epsilon
\ni
\Big(q-\sum^{N-1}_{j=0}q_j\Big)\odot p -r^{\odot N} 
=
q\odot p  - 1. 
\end{align*}

{\bf Right:} Put $r = 1- p \odot q_0 $. Define the symbol sequence $q_j = q_0 \odot r^{\odot j}$ for $j\in \mathbb{N}_0$. Then put $q  \sim \sum_{j=0}^\infty q_j $ with $q \in S^{-d}_\epsilon $, and write
\begin{align*}
S^{-N}_\epsilon
\ni
p\odot \Big(q-\sum^{N-1}_{j=0}q_j\Big) - r^{\odot N} 
=
p \odot q  - 1.
\end{align*}
\end{proof}


\begin{theorem} \label{thm:realansymbelliptic}
$p$ is elliptic in $S^d_\epsilon$ if and only if there is a finite $F\subset \widehat{G}$ so that:
\medskip
\begin{enumerate}

\item $p(z,\xi)$ is invertible for all $(z,[\xi])\in G_\epsilon \times (\widehat{G} \setminus F)$.

\medskip

\item The family of inverses satisfy
\begin{align*}
\sup_{(z,[\xi])\in G_\epsilon \times (\widehat{G}\setminus F )}\langle \xi \rangle^{d} ||p(z,\xi)^{-1} ||  < \infty.
\end{align*}

\end{enumerate}
\end{theorem}
\begin{proof}
Assume $p$ has a parametrix in $S^d_\epsilon$. That is, $q_0\in S^{-d}_\epsilon$ with $pq_0 - 1 \in S^{-1}_\epsilon$.
In that case, there is a $C>0$ so that, uniformly in $(z, [\xi]) \in G_\epsilon \times \widehat{G}$, we have
\begin{align*}
|| p(z,\xi) q_0(z,\xi) - I || \leq C \langle \xi \rangle^{-1},
\end{align*}
and so $(pq_0)(z,\xi)$ is invertible for $\langle \xi \rangle\geq R>0$. So is $p$, and
\begin{align*}
||p(z,\xi)^{-1}|| \leq ||q_0(z,\xi)|| \, \big|\big| (p(z,\xi)q_0 (z,\xi))^{-1} \big|\big|
\leq C' \langle \xi \rangle^{-d},
\end{align*}
where $C'>0$, and the inequality holds for $[\xi] \in \widehat{G}$ except the finite set with $\langle \xi \rangle < R$.
Conversely, if the two points above hold, consider the bounded subset $\{ p_Y \}_{Y \in \overline{\mathfrak{g}_\epsilon}}\subset S^d$.
It fulfils the hypotheses of 
Lemma~\ref{lmm:invmatsymb} with 
$K = \overline{\mathfrak{g}_\epsilon}$ and the above estimate over $z\in G_\epsilon$.
Thus, if we put $\chi_F(\xi) = 1_F([\xi])I_{d_\xi}$ for all $[\xi]\in \widehat{G}$, then $\chi_F \in S^{-\infty}_\epsilon$, it gives
\begin{align*}
(\chi_F + (1-\chi_F) p)^{-1} \in S^{-d}_\epsilon,
\end{align*}
and so $p\in S^d_\epsilon$ has a parametrix in $S^d_\epsilon$.
\end{proof}

\newpage
As a consequence, elliptic left-invariant operators admit a parametrix in $\Op \, S_\epsilon$. More generally, only the leading term in $p$ has to belong to such an operator:

\begin{proposition}
Suppose that $p$ of $P= \Op(p)\in \Psi^k(G)$ is of the form
\begin{align*}
p = p_k  + p_{k-1},
\end{align*}
where the following holds:
\medskip
\begin{enumerate}

\item $p_k \in S^k$ is elliptic and independent of $x\in G$.

\medskip

\item $p_{k-1} \in S^{k-1}_\epsilon$ for a fixed $\epsilon>0$.

\end{enumerate}
\medskip
Then $P$ has a parametrix in $\Op \, S^{-k}_\epsilon$ for this particular $\epsilon>0$.
\end{proposition}
\begin{proof}
Observe that 
\begin{align*}
\sup_{(x,Y)\in G \times \mathfrak{g}_\epsilon} \langle \xi \rangle^{-k} ||p_Y (x,\xi) -  p_k(\xi)|| = O_{\lambda_\xi \to \infty}( \langle \xi \rangle^{-1} ),
\end{align*}
and the symbols $p_Y(x,\xi)$ are invertible for all $(x,Y)\in G \times \mathfrak{g}_\epsilon$ if $\langle \xi \rangle$ is large enough. 
This is because $p_k(\xi)$ is invertible for all sufficiently large $\langle \xi \rangle$
by \cite[Theorem 4.1]{RuzhanskyTurunenWirth2014}.
Furthermore, it gives a finite $F \subset \widehat{G}$ and a $C>0$ such that 
\begin{align*}
|| \langle \xi \rangle^{k} p_k (\xi)^{-1} || \leq C
\quad
\textnormal{for all}
\quad
[\xi]\in \widehat{G} \setminus F,
\end{align*}
and for $\langle \xi \rangle$ suitably large, we can then write
\begin{align*}
\langle \xi \rangle^{k}  ||p_Y (x,\xi)^{-1}||
&=
\Big|\Big| \Big( I + p_k (\xi)^{-1} p_{k-1}(\exp(iY)x,\xi)  \Big)^{-1} \langle \xi \rangle^{k} p_k (\xi)^{-1}    \Big|\Big| \\
&\leq
\Big( 1 - || p_k (\xi)^{-1} p_{k-1}(\exp(iY)x,\xi)||  \Big)^{-1} || \langle \xi \rangle^{k} p_k (\xi)^{-1} || \\
&\leq
\Big( 1 -  C \langle \xi \rangle^{-k}  || p_{k-1}(\exp(iY)x, \xi)||  \Big)^{-1} C,
\end{align*}
where the inner term is smaller than $\frac{1}{2}$ for all $(x,Y)\in G \times \mathfrak{g}_\epsilon$ if $\langle \xi \rangle$ is large enough.
Therefore the conditions of Theorem~\ref{thm:realansymbelliptic} are fulfilled.
\end{proof}

In other words, the "leading term"  determines if $P$ has the property that we seek. These observations indicate that $\Op\, S^d_\epsilon$ contain many non-trivial operators.

\begin{example}
{\rm 
Operators in $\Psi(\mathbb{T}^n)$, acting on $u\in C^\infty(\mathbb{T}^n)$, are of the form
\begin{align*}
\Op(p) u(x)
= \sum_{k\in \mathbb{Z}^n} e^{ik\cdot x} p(x,k) \int_{\mathbb{T}^n} u(y) e^{-ik\cdot y} \, dy, 
\end{align*}
and we may take $p$ to be the symbol 
\begin{align*}
p(x,k) = |k|^2 + \sum_{j=1}^n g_j(x) k_j + f(x),
\end{align*}
where $\{ g_j \}_{j=1}^n$, $f$ are bounded holomorphic on $\mathbb{T}^n_\epsilon$ ($2\pi$-periodic on the polystrip $\mathbb{R}^n_\epsilon$).
Quantization of this particular symbol gives the Laplacian $\Delta_{\mathbb{T}^n}$ plus lower order terms. The result above says that $p$ is elliptic in $ S^2_\epsilon (\mathbb{T}^n \times \mathbb{Z}^n)$, e.g. via $p(x,k)^{-1}$ for large $k$, which is not a sum of products of functions in $x$ and $k$ separately.
}
\end{example}

\newpage
\section{Outlook} 
The Boutet de Monvel theorem/conjecture is yet not fully proven, but it has important implications for differential operators on real-analytic manifolds.
We expect that the results on Lie groups extend to compact homogeneous spaces, where the chosen metric is related to a bi-invariant metric on the group acting on it.
Our hope is that an analogous conjecture can be made for manifolds with boundary, and, if true, that it reveals a path to studying real-analytic boundary value problems in Sobolev spaces of holomorphic functions.

\section{Acknowledgements}
The author wishes to thank the anonymous reviewer for pointing out a few mistakes in the original manuscript.

\bibliographystyle{siamplain}
\bibliography{references}
\end{document}